\documentclass[]{article}
\usepackage{fullpage}
\usepackage{amssymb,amsmath}
\usepackage{subcaption}
\usepackage{graphicx}
\usepackage{longtable,booktabs}
\usepackage{multirow}
\usepackage{lipsum}
\usepackage{amsfonts}
\usepackage{graphicx}
\usepackage{epstopdf}
\usepackage{algorithmic}
\usepackage{longtable,booktabs}
\usepackage{multirow}
\usepackage{color}
\usepackage{caption}
\usepackage{subcaption}

\usepackage{cleveref}
\usepackage{algorithm}
\usepackage{algorithmic}
\usepackage{url}

\newcommand{\dee}[2]{\frac{\textrm{d} #1}{\textrm{d} #2}}
\newcommand{\YA}[1]{\textcolor{black}{#1}}

\title{Cluster Gauss-Newton method for finding multiple approximate minimisers of nonlinear least squares problems with applications to parameter estimation of pharmacokinetic models}
  
\author{Yasunori Aoki\thanks{Visiting Associate Professor, National Institute of Informatics, 2-1-2 Hitotsubashi, Chiyoda, Tokyo 101-8430, Japan, and Visiting Researcher, Sugiyama Laboratory, RIKEN Baton Zone Program, 1-7-22 Suehiro-cho, Tsurumi-ku, Yokohama, Kanagawa 230-0045, Japan (yaoki@uwaterloo.ca).} \and Ken Hayami\thanks{National Institute of Informatics, and The Graduate University for Advanced Studies (SOKENDAI ), 2-1-2 Hitotsubashi, Chiyoda, Tokyo 101-8430, Japan (hayami@nii.ac.jp). The author was supported in part by JSPS KAKENHI Grant Number 15K04768 and 15H02968.} \and Kota Toshimoto\thanks{Sugiyama Laboratory, RIKEN Baton Zone Program, 1-7-22 Suehiro-cho, Tsurumi-ku, Yokohama, Kanagawa 230-0045, Japan.} \and Yuichi Sugiyama\thanks{Sugiyama Laboratory, RIKEN Baton Zone Program, 1-7-22 Suehiro-cho, Tsurumi-ku, Yokohama, Kanagawa 230-0045, Japan.}}

\begin{document}
\maketitle
%
%

\begin{abstract}
	Parameter estimation problems of mathematical models can often be formulated as nonlinear least squares problems.  Typically these problems are solved numerically using iterative methods.  The \YA{local minimiser} obtained using these iterative methods usually depends on the choice of the initial iterate. Thus, the estimated parameter and subsequent analyses using it depend on the choice of the initial iterate.  One way to reduce the analysis bias due to the choice of the initial iterate is to repeat the algorithm from multiple initial iterates (i.e. use a multi-start method).  However, the procedure can be computationally intensive and is not always used in practice.  To overcome this problem, we propose the Cluster Gauss-Newton (CGN) method, an efficient algorithm for finding multiple \YA{approximate minimisers} of nonlinear-least squares problems.  \YA{CGN} simultaneously solves the nonlinear least squares problem from multiple initial iterates.  \YA{Then, CGN} iteratively improves the solutions from these initial iterates similarly to the Gauss-Newton method. However, it uses a global linear approximation instead of the Jacobian.  The global linear approximations are computed collectively among all the iterates to minimise the computational cost. 
We use physiologically based pharmacokinetic (PBPK) models used in pharmaceutical drug development to demonstrate its use and show that CGN is computationally more efficient and more robust against local minima compared to the standard Levenberg-Marquardt method, as well as state-of-the art multi-start and derivative-free methods.
\end{abstract}

\section{Introduction}\label{header-c1607}

The parameter estimation of mathematical models often boils down to solving nonlinear least squares problems. Hence, algorithms for solving nonlinear least squares problems are widely used in many scientific fields. 

The most traditional least squares solver is the Gauss-Newton method~\cite{gauss1857theory, Bjorck1996}.  In practice, the Gauss-Newton method with regularisation (i.e., Levenberg-Marquardt (LM) method~\cite{marquardt1963algorithm, more1978levenberg}) or with the Trust-Region method~\cite{conn2000trust} is often used.
Recently, derivative-free methods, which do not explicitly use derivative information of the nonlinear function, have been developed.  These methods are usually computationally more efficient as it avoids the costly computation of the derivatives. Also, they can be applied even to  problems where the mathematical models are `black box'.  The state of the art derivative-free algorithms are DFO-LS~\cite{cartis2017derivative} and POUNDERS~\cite{wild2017chapter}.  A comprehensive review of the derivative-free methods can be found in \cite{larson2019derivative}.

Another approach for obtaining a solution of nonlinear least squares problems is to directly minimise the sum of squared residuals (SSR) using generic optimisation algorithms for the scalar objective function.  The most classical approach is to obtain a solution where the gradient of the SSR becomes zero using the Newton method.  As it is usually too costly to compute the Hessian of the SSR, Quasi-Newton methods which approximate the Hessian are used.  The commonly used Quasi-Newton method is the BFGS method~\cite{broyden1970convergence, fletcher1970new, goldfarb1970family, shanno1970conditioning, shanno1970optimal}.  Another approach which makes use of the Newton-type method for optimisation is Implicit Filtering \cite{kelley2011implicit} which combines grid search and the Newton method.  In addition to these optimisation algorithms, we can use numerous global optimisation algorithms when bound constraints are given.  For example, Surrogate Optimisation \cite{gutmann2001radial}, Genetic Algorithm \cite{goldberg1988genetic}, Particle Swarm algorithm \cite{kennedy1995particle, mezura2011constraint}, and DIRECT \cite{jones1993lipschitzian} are well known global optimisation algorithms. 

Although there are a variety of algorithms to solve the nonlinear least squares problems as listed above, they mostly focus on finding one solution. To the best of our knowledge, there is very limited methodological development on algorithms for simultaneously finding multiple \YA{approximate minimiser}s of nonlinear least squares problems.
 \YA{For instance,} when using the Levenberg-Marquardt method, the local algorithm often gives a local minimiser which depends on the choice of the initial iterate. To reduce the analysis bias due to the initial iterate used for a local algorithm, it is a good practice to repeatedly use the local algorithm with various initial iterates, as in multi-start methods~\cite{boender1982stochastic}.
Similarly, for problems where bound constraints are given, one can use global optimisation algorithms to find one of the global minimisers.  On the other hand, if there are multiple global minisers, the global minimiser found can depend on the algorithm setting, for example, the random seed.  Hence, it is beneficial to use global optimisation algorithms with various settings repeatedly if the uniqueness of the global minimiser is not guaranteed.
 \YA{The trivial bottleneck of repeatedly using these algorithms is the computation cost. 
 In this paper, we propose a new method addressing this computational challenge of finding multiple \YA{approximate minimisers
 of } nonlinear least squares problems.}

Our algorithm development for finding multiple \YA{local minimiser}s of nonlinear least squares problems was motivated by a mathematical model of pharmaceutical drug concentration in a human body called the physiologically based pharmacokinetic (PBPK) model \YA{\cite{Watanabe2009}}.  The PBPK model is typically a system of mildly nonlinear stiff ordinary differential equations (ODEs) with many parameters. This type of mathematical model is constructed based on the knowledge of the mechanism of how the drug is absorbed, distributed, metabolised and excreted. Given the complexity of this process and the limitation of the observations we can obtain from a live human subject, the model parameters cannot be uniquely identified from the observations, meaning that there are non-unique global minimisers to the nonlinear least squares problem. The estimated parameters of the PBPK model are used to simulate the drug concentration of the patient from whom we are often unable to test the drug on (e.g., children, pregnant person, a person with rare genetic anomaly) or to predict the experiment that is yet to be run (different amount of drug administration, multiple drug used at the same time). As the simulated drug concentration is used to predict the safety of the drug in these different scenarios, it is essential to consider multiple predictions based on multiple possible parameters that are estimated from the available observations.  A motivating example is presented in Appendix~\ref{sec::motivatingExample}. \YA{Another reason why we want to obtain multiple sets of parameters is that we can understand which parameter cannot be estimated from the available data.  This will motivate the pharmaceutical scientists to perform additional (e.g., in-vitro or in-animal) experiments to determine these parameters that were not estimable from the available data.}

In \cite{Aoki2011, Aoki2014} we proposed the Cluster Newton (CN) method for obtaining multiple solutions of {\bf a system of \YA{underdetermined} nonlinear equations}.
In recent years CN has been used in the field of pharmaceutical science \cite{Yoshida2013, Fukuchi2017, Asami2017, Toshimoto2017, Kim2017, Nakamura2018}. 
For example, \cite{Toshimoto2017} used the parameters estimated by CN to predict the adverse drug effect, \cite{Nakamura2018} used the estimated parameters to predict the outcome of a clinical trial.  
However, based on these applications of CN, we observed the necessity to formulate the problem as a nonlinear least squares problem and make the method more robust against noise in the observed data. This is mainly because actual pharmaceutical data may contain measurement error and inconsistency coming from an inadequate model.  (See Appendix~\ref{header-c1650}.)

\subsection{Nonlinear least squares problem of our interest}\label{header-c1659}

In this paper, we propose an algorithm for obtaining multiple \YA{approximate minimiser}s of nonlinear least squares problems 
\begin{eqnarray}\min_{\boldsymbol{x}}|| \boldsymbol{f}(\boldsymbol x)-\boldsymbol{y}^*||_2^{\,2}
\label{eq:mathProblem}
\end{eqnarray}
which do not have a unique \YA{solution (global minimiser)}, that is to say, there exist $\boldsymbol x^{(1)}\neq\boldsymbol x^{(2)}$ such that
\begin{eqnarray}
\min_{\boldsymbol{x}}|| \boldsymbol{f}(\boldsymbol x)-\boldsymbol{y}^*||_2^{\,2}=|| \boldsymbol{f}(\boldsymbol x^{(1)})-\boldsymbol{y}^*||_2^{\,2}=|| \boldsymbol{f}(\boldsymbol x^{(2)})-\boldsymbol{y}^*||_2^{\,2} \,.
\end{eqnarray}

Here, $\boldsymbol{f}$ is a nonlinear function from $\mathbb{R}^n$ to $\mathbb{R}^m$, $\boldsymbol x^{(1)},\,\boldsymbol x^{(2)}\in \mathbb{R}^n$ and $\boldsymbol{y}^*\in \mathbb{R}^m$. The nonlinear function $\boldsymbol{f}$ can be derived from a mathematical model, the vector $\boldsymbol x$ can be regarded as a set of model parameters which one wishes to estimate, and the vector $\boldsymbol{y}^*$ can be regarded as a set of observations one wishes to fit the model to.
\YA{In the motivating problem of estimating parameters of PBPK models, since there are usually insufficient observations, the corresponding nonlinear least squares problem has multiple global minimisers.  Therefore, we are interested in finding multiple minimisers instead of just one minimiser.} 


We shall call the following quantity the sum of squared residuals (SSR): 
\begin{eqnarray}
|| \boldsymbol{f}(\boldsymbol x)-\boldsymbol{y}^*||_2^{\,2}\,, \label{eq::SSR}
\end{eqnarray}
and use it for the quantification of the goodness of $\boldsymbol x$ as the approximation of the \YA{solution} of the least squares problem \eqref{eq:mathProblem}.

\subsection{Well known examples in pharmacokinetics}
\begin{figure}[h]
\begin{subfigure}[b]{0.49\columnwidth}
\includegraphics[width=55mm]{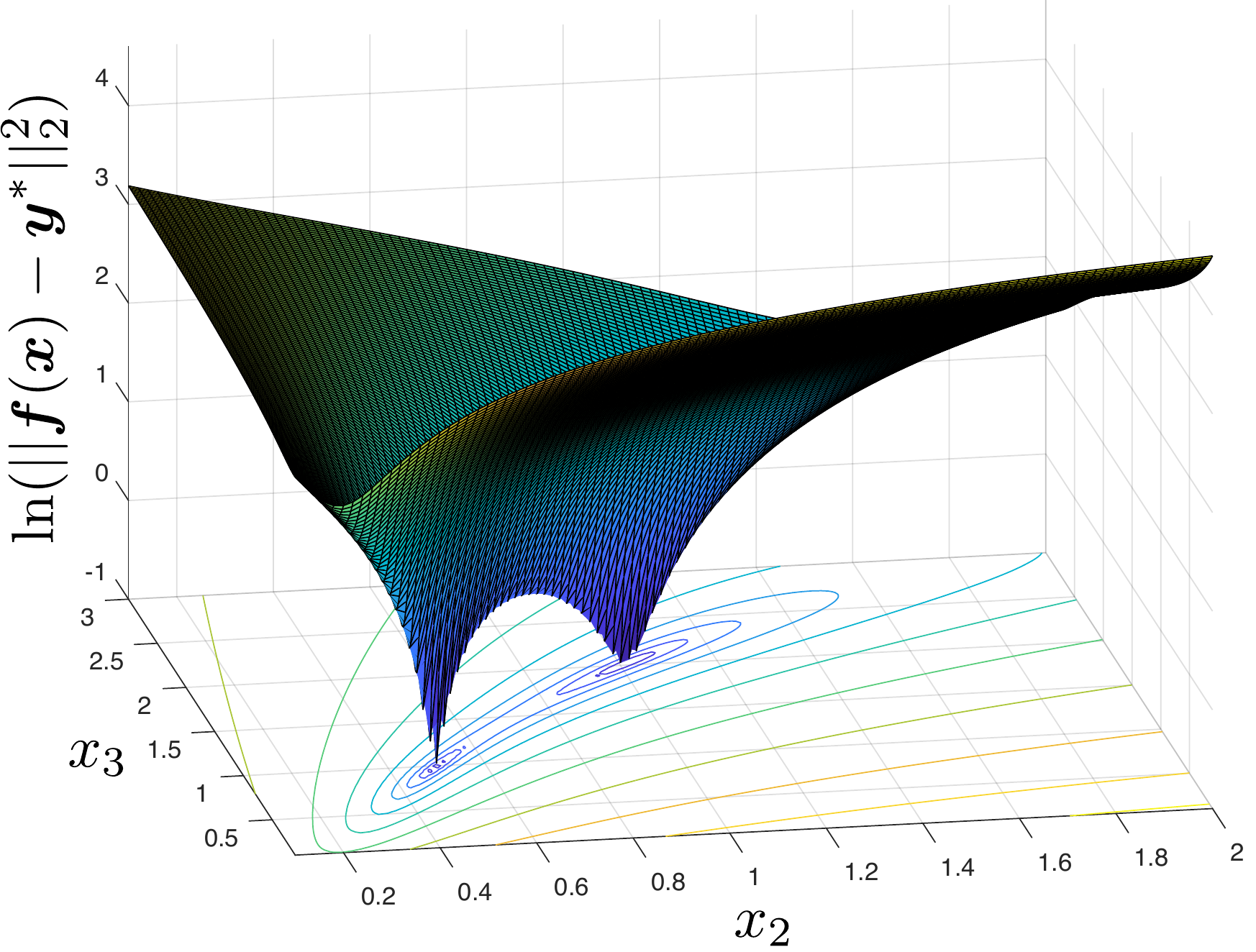}
\caption{Flip-flop kinetics.}
\label{Fig::surfPlot_flipflop}
\end{subfigure}
\begin{subfigure}[b]{0.49\columnwidth}
\includegraphics[width=55mm]{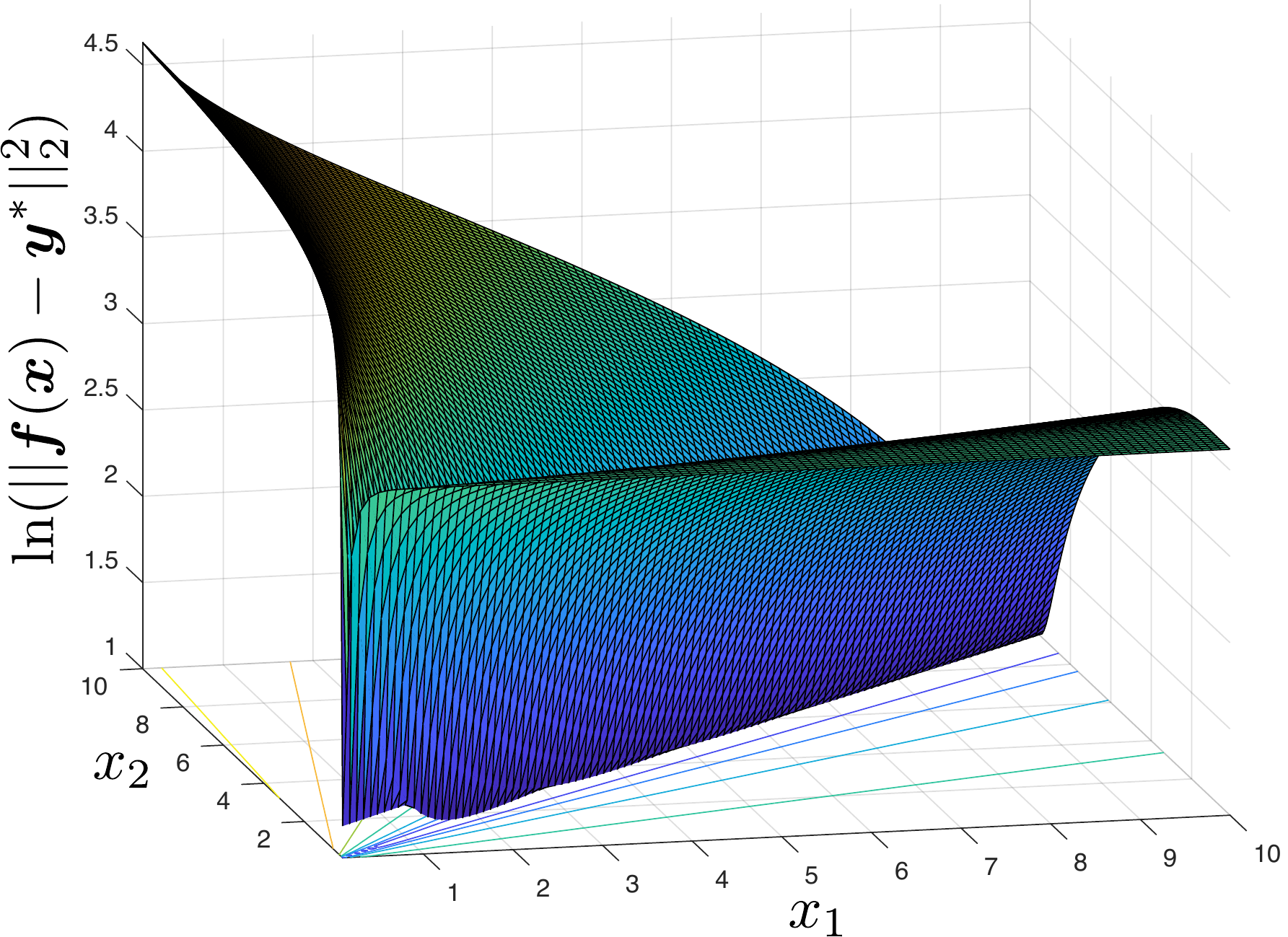}
\caption{Total amount of the drug in the body as observation.}
\label{Fig::surfPlot_nonidentifiable}
\end{subfigure}
\caption{Surface plots of the sum of squared residuals for pharmacokinetics model parameter estimation problems with non-unique global minimisers.}
\label{Fig::surfacePlots}
\end{figure}

\YA{In this subsection, we present two simple pharmacokinetics parameter estimation problems where the corresponding nonlinear least squares problems have non-unique global minimisers.  
The first example is called `flip-flop kinetics' and can be found in most   standard textbooks in pharmacokinetics (e,g, \cite{gibaldi1982drugs}).  Flip-flop kinetics occurs when estimating the pharmacokinetic parameters of the drug that is orally given (e.g., as a pill or a tablet) to patients based on the observation of the drug concentration in the blood plasma.  The simplest mathematical model for this pharmacokinetics can be written as follows:}
\begin{eqnarray}
\dee{u_1}{t}&=&-Ka\,u_1\\
\dee{u_2}{t}&=&\frac{Ka\,u_1-CL\,u_2}{V}\\
u_1(t=0)&=&100\\
u_2(t=0)&=&0
\end{eqnarray}
where
\begin{eqnarray}
	CL&=&10^{x_1}\\
Ka&=&10^{x_2}\\
V&=&10^{x_3}\,.
\end{eqnarray}
\YA{For this problem $u_2$ corresponds to the drug concentration in the blood plasma, which is observable. It can be shown analytically that there are two distinct parameter sets that realise the same drug concentration time-course curve.  Figure~\ref{Fig::surfPlot_flipflop} shows the surface plot of the sum of squared residuals of the corresponding nonlinear least squares problem.  As can be seen, there are two global minimisers for this nonlinear least squares problem.}

\YA{The second example is when an isotope labelled substance is given to patients via intervenes administration and pharmacokinetics parameters are estimated from the total amount of the substance in the body (i.e., not from the concentration).  The simplest mathematical model for this pharmacokinetics can be written as follows:}
\begin{eqnarray}
\dee{u}{t}&=&-\frac{CL}{V}u\\
u(t=0)&=&100
\end{eqnarray}
where
\begin{eqnarray}
	CL&=&10^{x_1}\\
	V&=&10^{x_2}.
\end{eqnarray}

For this problem, $u$ corresponds to the total amount of the substance in the body, which is the observable quantity. In this case, there are infinite sets of parameters $(x_1, x_2)$ which satisfy a linear relation, that realises the given time-course substance amount observation. Figure~\ref{Fig::surfPlot_nonidentifiable} is a surface plot of the sum of squared residuals of the corresponding nonlinear least squares problem. As can be seen, there is a line of global minimisers for this nonlinear least squares problem.

\YA{In this subsection, we have shown the simplest possible forms of this issue of non-unique global minimisers in pharmacokinetics.  Thus, we would like to point out that one cannot assume the uniqueness of the global minimiser for more complex pharmacokinetic models, for example, the PBPK model of our interest.}

\section{Method: Algorithm}\label{header-c1683}
In this section we describe the proposed algorithm.  We first introduce a rough concept using a toy example in Subsection~\ref{sec:briefAlg} and then introduce the full algorithm in detail in Subsection~\ref{sec:detailAlg}.

\subsection{Brief explanation of the algorithm}\label{sec:briefAlg}
The aim of the proposed Cluster Gauss-Newton (CGN) algorithm is to efficiently find multiple \YA{approximate minimiser}s of the nonlinear least squares problem \eqref{eq:mathProblem}.  We do so by first creating a collection of initial guesses which we call the `cluster'.  Then, we move the cluster iteratively using linear approximations of the nonlinear function $\boldsymbol{f}$, similarly to the Gauss-Newton method \cite{Bjorck1996}.

The unique idea in the CGN method is that the linear approximation is constructed collectively throughout the points in the cluster instead of using the Jacobian matrix which approximates the nonlinear function linearly at a point as in the Gauss-Newton or LM. By using points in the cluster to construct a linear approximation, instead of explicitly approximating the Jacobian, we minimise the computational cost for each iteration. \YA{In addition, by constructing linear approximation using non-local information, CGN is more likely to converge to approximate minimisers with smaller SSR than methods using the Jacobian.}

\begin{figure}
\includegraphics[width=80mm]{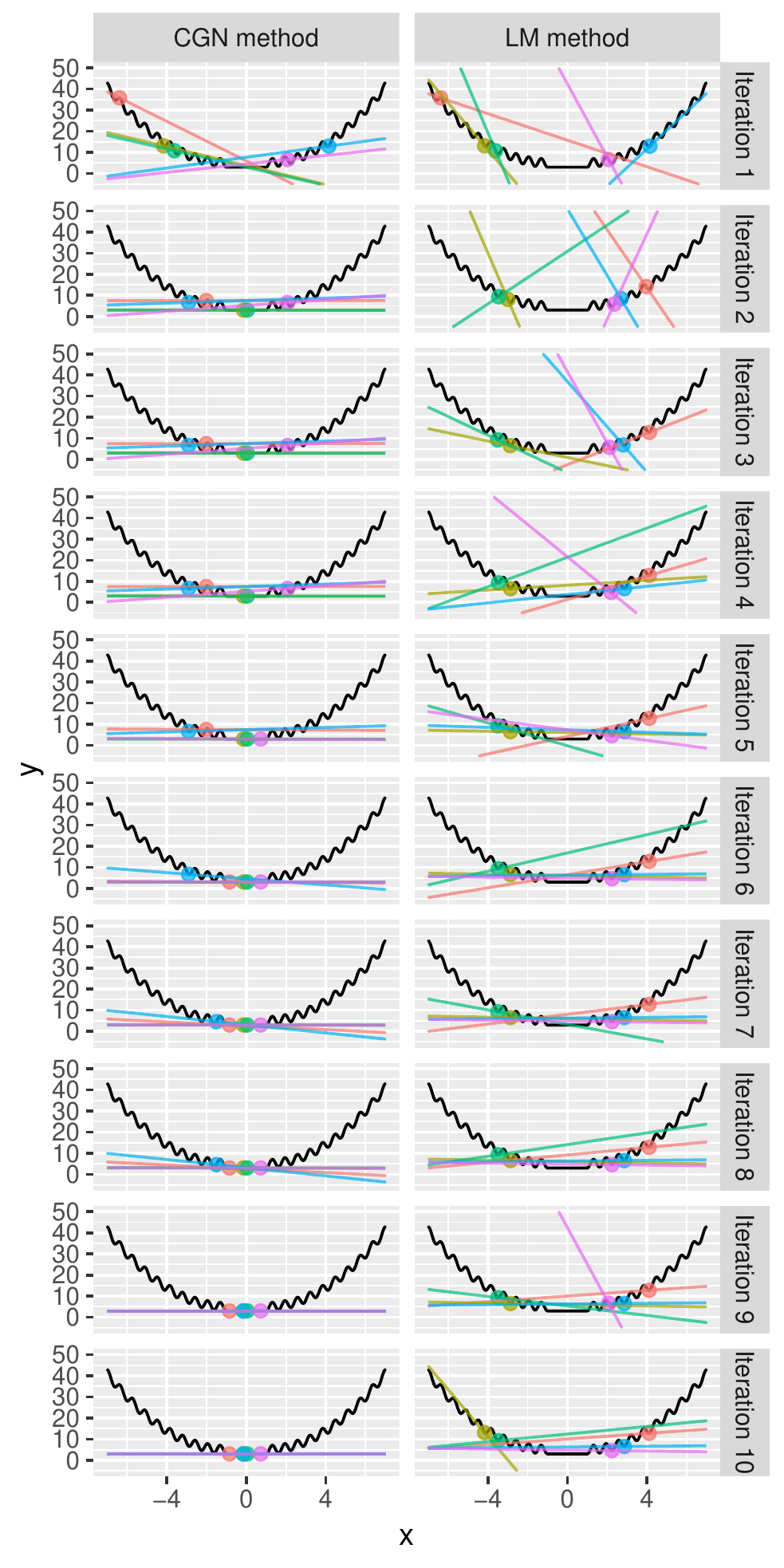}
\caption{Schematic comparison of CGN and LM.  Dots represent the iterates with their function values in each iteration and the lines represent linear approximations used to update the iterates.  As can be seen in this figure, the linear approximation used in CGN follows the global trend of the function, while the slope used in LM captures the local feature of the function.  As a result, when using CGN all the iterates converge to the minimisers with the smallest residual (in this case global minimisers), while all the iterates of LM converge to local minimisere whose residuals are not the minimum. In addition, for the 10 iterations presented in the figure, CGN required only 50 function evaluations while LM required 121 function evaluations since the slope is approximated using finite differences.}
\label{Fig::drawingComparisonLA}
\end{figure}

In order to visualise the key differences between the proposed linear approximation (CGN) and the Jacobian (i.e., derivative) approaches, we consider the  nonlinear function 
\begin{eqnarray}
	f&=&\left\{\begin{array}{ll}
       (x+1)^2-2\cos(10(x+1))+5&\textrm{if }x<-1\\
       3&\textrm{if }-1 \leq x \leq 1 \\
	(x-1)^2-2\cos(10(x-1))+5&\textrm{if }x>1	
\end{array}	\right.
\end{eqnarray}
 (see Figure~\ref{Fig::drawingComparisonLA}) and aim to find global minimisers. Any point $x\in[-1,1]$ is a global minimiser of this problem. Hence, this problem has nonunique global minimisers.  Let the points of the initial iterates be:
 
\begin{eqnarray}
&&x_1=-6.3797853, \qquad x_2=-4.1656025, \qquad x_3=-3.6145728,\nonumber\\ &&x_4=2.0755468, \qquad x_5= 4.1540421.
\end{eqnarray}
We now compute the linear approximations used to move these points in the cluster to minimise the function $f$.

\noindent\underline{Gradient (LM)}

For this nonlinear function, \YA{since the function is given in analytic form,} we can obtain the gradient explicitly.  \YA{In,} practice, when $f$ is given as a ``black box'', we can approximate the derivative by a finite difference scheme, for example, $f'(x_i)\approx\frac{f(x_i+\epsilon)-f(x_i)}{\epsilon}$.  Then, the linear approximation at $x_i$ can be written as \\
$f(x) \approx \frac{f(x_i+\epsilon)-f(x_i)}{\epsilon}(x-x_i)+f(x_i)$. Notice that it requires one extra evaluation of $f$ at $x_i+\epsilon$ for each $x_i$. This number of extra function evaluation is, when evaluating a full gradient estimate, equal to the number of independent variables of $f$.  (This is not the case if a directional derivative estimate is used.)
  \YA{If $f$ is given by a system of ODEs, one may use the
adjoint method to obtain the derivatives more efficiently. However, it requires solving an additional system of ODEs (the adjoint equation).}
More importantly, iterates of methods based on the gradient may converge to local minimisers, since they use local gradient information.

\noindent\underline{Cluster Gauss-Newton (CGN) method (proposed method)}

In the proposed method, we construct a linear approximation for each point in the cluster while using the value of $f$ at other points in the cluster to globally approximate the nonlinear function with a linear function.  The influence of another point in the cluster to the linear approximation is weighted according to how close the point is to the point of approximation, i.e., 
\begin{eqnarray}
	\min_{a_{(i)}}\sum_{i\neq j}\left(d_{j(i)}\left((x_j-x_i)a_{(i)}-\left(f(x_j)-f(x_i)\right)\right)\right)^{\,2}\label{eq:linApproxNCN}
\end{eqnarray}
where $a_{(i)}$ is the slope of the linear approximation at $x_i$ and the linear approximation at $x_i$ can be written as $f(x)\approx a_{(i)}(x-x_i)+f(x_i)$.  There are many possibilities for the weight function $d_{j(i)}$.  In this paper, we choose $d_{j(i)}=(x_j-x_i)^{-2\gamma}$ where $\gamma \geq 0$ ($\gamma = 0$ corresponds to uniform weight).
Note that Equation\eqref{eq:linApproxNCN} can also be regarded as a weighted least squares solution of a system of linear equations \YA{where the weight is chosen to be $d_{j(i)}$}.  The weight is motivated by the fact that we weight the information from the neighbouring points in the cluster more than the ones further away when constructing the linear approximation. Note that we do not require any extra evaluation of $f$ for obtaining these linear approximations.

For the multi-dimensional nonlinear function $\boldsymbol{f}:\mathbb{R}^n\rightarrow \mathbb{R}^m$, when we have $N$ points in the cluster, we solve the following linear least squares problem:
\begin{eqnarray}
\min_{A_{(i)}\in \mathbb{R}^{m\times n}} \left|\left|D_{(i)}\left(\Delta X_{(i)} A^\textrm{T}_{(i)} -\Delta Y_{(i)}\right) \right|\right|_\textrm{F}	\label{eq::linApprox_multiD_rough}
\end{eqnarray}
where $D_{(i)}=\textrm{diag}(d_{1i},...,d_{Ni})$ is a diagonal matrix defining the weights, $\Delta X_{(i)}\in\mathbb{R}^{N\times n}$ is the difference between all the cluster points and $\boldsymbol{x}_i$, and $\Delta Y_{(i)}\in\mathbb{R}^{N\times m}$ is the difference between the nonlinear function $\boldsymbol{f}$ evaluated at all the cluster points and at $\boldsymbol{x}_i$. The precise definition of these quantities and derivation of \eqref{eq::linApprox_multiD_rough} is given in the next subsection.
 
For the one dimensional case (i.e., $m=n=1$), the linear approximations at each point for the first ten iterations for both CGN and LM are shown in Figure~\ref{Fig::drawingComparisonLA}.  As can be seen, the gradient used in LM captures the local behaviour of the nonlinear function. The linear approximation used in CGN, on the other hand, captures the global behaviour of the nonlinear function.  After nine iterations of the CGN, all the points reached the minimisers with smallest SSR.  On the other hand, the LM converges to local minimisers whose SSR are not necessarily the minimum.

\subsection{Detailed description of the algorithm}\label{sec:detailAlg}
Next, we describe the proposed CGN algorithm in detail. 
In this subsection, we denote a scalar quantity by a lower case letter e.g., $a$, $c$, a matrix by a capital letter, e.g., $A$, or $M$, and a column vector by a bold symbol of a lower case letter, e.g., $\boldsymbol{v}$, $\boldsymbol{a}$, unless otherwise specifically stated.  Super script $\textrm{T}$ indicates the transpose. Hence, $\boldsymbol{v}^{\textrm{T}}$ and $\boldsymbol{a}^{\textrm{T}}$ are row vectors.\\
\\
{\setlength{\parindent}{0.5cm}

\setlength{\leftskip}{0.5cm}

\noindent{\bf 1) Pre-iteration process}\label{header-c1684}

{\bf 1-1) Create initial cluster}

 The initial iterates of CGN, a set of vectors
$\{\boldsymbol{x}_{i}^{(0)}\}_{i=1}^{N}$ are genrated using uniform random numbers in each component 
 within the domain of initial estimate of the plausible location of \YA{global minimisers} given by the user. The unique point of the CGN is that the user specifies the domain of the initial estimates instead of a point.  In this paper, we assume that the domain of initial guess is given by the user by two sets of vectors 
$\boldsymbol{x}^\textrm{L}$, $\boldsymbol{x}^\textrm{U}$, and the value
$x_{ji}^{(0)}$ is sampled from the uniform distribution between $x^\textrm{L}_j$ and $x^\textrm{U}_j$, where $x_{ji}^{(0)}$ is the $j$ th element of vector $\boldsymbol{x}_i^{(0)}$.

Note that this does {\bf not} mean that all the following iterates $\boldsymbol{x}_i^{(k)} \hspace{2mm} (k \ge 1, i=1,2, \ldots, N)$ must satisfy 
$ \boldsymbol{x}^L \le \boldsymbol{x}_i^{(k)} \le \boldsymbol{x}^U$.

\setlength{\leftskip}{1cm}
Store the initial set of vectors in a matrix $X^{(0)}$, i.e.,

\begin{eqnarray}X^{(0)} =[\boldsymbol{x}_{1}^{(0)},
\boldsymbol{x}_{2}^{(0)}, 
\ldots ,
\boldsymbol{x}_{N}^{(0)}
]\label{eq:initMatrixX}\end{eqnarray}

\noindent where the super script (0) indicates the initial iterate.

Evaluate the nonlinear function $\boldsymbol{f}$ at each
$\boldsymbol{x}_{i}^{(0)}$ as $\boldsymbol{y}_{i}^{(0)}=\boldsymbol{f}(\boldsymbol{x}_{i}^{(0)}) \hspace{2mm} (i=1,2,\ldots,N)$ and store in matrix $Y^{(0)}$, i.e.,

\begin{eqnarray}Y^{(0)}=\left[ \boldsymbol{y}_{1}^{(0)}, \boldsymbol{y}_{2}^{(0)}, \ldots, \boldsymbol{y}_{N}^{(0)}\right]. \label{eq:initMatrixY}\end{eqnarray}

If the function $\boldsymbol{f}$ cannot be evaluated at
$\boldsymbol{x}_{i}^{(0)}$, then re-sample
$\boldsymbol{x}_{i}^{(0)}$ until $\boldsymbol{f}$ can be evaluated.

Compute the sum of squared residuals vector $\boldsymbol{r}^{(0)}$
i.e.,

\begin{eqnarray}
r_{i}^{0}=||\boldsymbol{y}_i^{0} - \boldsymbol{y}^* ||_2^2 \hspace{3mm} (i=1,2,\ldots,N)\\
\boldsymbol{r}^{(0)}=\left[ r_{1}^{(0)}, r_{2}^{0}, \ldots, r_{N}^{0} \right]^{\textrm{T}}. \label{eq:initResidualVector}\end{eqnarray}

  \YA{The concise pseudo-code for the creation of the initial cluster can be found in Algorithm~\ref{alg:initialCluster}.}
  
\begin{algorithm}     
                   
\caption{$(X^{(0)},Y^{(0)},\boldsymbol{r}^{(0)})$={\bf create\_initialCluster}$(\boldsymbol{x}^\textrm{L},\boldsymbol{x}^\textrm{U}, N, \boldsymbol{f},\boldsymbol{y}^*)$}         
\label{alg:initialCluster}                          
\begin{algorithmic}                  
\STATE{$X^{(0)}\gets$ allocate memory for $n \times N$ matrix}
\STATE{$Y^{(0)}\gets m \times N$ matrix of not a number (NaN)}

\FORALL{$i= 1,...,N$} 
 
 \WHILE{$\boldsymbol{y}^{(0)}_i$ is not a NaN vector}
 \FOR{$j\gets 1, n$} 
 \STATE {Draw random number $x^{(0)}_{ij}$ from $\mathcal{U}(x^\textrm{L}_j, x^\textrm{U}_j)$} 
 \ENDFOR
 
  \IF{$\boldsymbol{f}(\boldsymbol{x}_i)$ can be evaluated}
 \STATE{
    $\boldsymbol{y}^{(0)}_i \gets \boldsymbol{f}(\boldsymbol{x}_i)$
  }
  \ELSE
  \STATE{
    $\boldsymbol{y}^{(0)}_i \gets $ NaN vector of length $m$
  }
 \ENDIF
 \ENDWHILE
 \STATE{
    $r_i^{(0)} \gets || \boldsymbol{y}_i^{(0)}-\boldsymbol{y}^*||_2^{\,2}$
  }
\ENDFOR
\end{algorithmic}
\end{algorithm}

{\bf 1-2) Initialise regularisation parameter vector}

Fill the regularisation parameter vector $\boldsymbol{\lambda}^{(0)}\in \mathbb{R}^N$, with
the user-specified initial regularisation parameter
$\lambda_\textrm{init}$ .i.e.,

\begin{eqnarray}
\boldsymbol{\lambda}^{(0)}=[\lambda_\textrm{init},\lambda_\textrm{init},...,\lambda_\textrm{init}]^{\textrm{T}}
\label{lambda_init}
\end{eqnarray}

\setlength{\leftskip}{0.5cm}
\noindent\textbf{2) Main iteration}

Repeat the following procedure until the user specified stopping
criteria are met. We denote the iteration number as $k$, which starts
from $0$ and is incremented by $1$ after each iteration.

\setlength{\leftskip}{1cm}

\textbf{2-1) Construct weighted linear approximations of the nonlinear
function}\\
We first construct a linear approximation around the point $\boldsymbol{x}_{i}^{(k)}$, s.t., 
\begin{eqnarray}\boldsymbol{f}(\boldsymbol{x})\approx A_{(i)}^{(k)}(\boldsymbol{x}-\boldsymbol{x}_{i}^{(k)})+\boldsymbol{f}(\boldsymbol{x}_{i}^{(k)}) \,,\end{eqnarray}

Here, $A^{(k)}_{(i)} \in \mathbb{R}^{m\times n}$ describes the slope of the linear approximation around $\boldsymbol{x}_{i}^{(k)}$.

The key difference of our algorithm compared to others is that we
construct a Jacobian like matrix $A_{(i)}^{(k)}$ collectively using all
the function evaluations of $\boldsymbol{f}$ in the previous
iteration, i.e., we solve
\begin{eqnarray}
\qquad \qquad \min_{A \in \mathbb{R}^{m \times n} }\sum_{j=1}^N\left[ d_{j(i)}^{(k)}   \left|\left| \boldsymbol{f}(\boldsymbol{x}_j^{(k)})-\left\{ A \left(\boldsymbol{x}_j^{(k)}-\boldsymbol{x}_i^{(k)}\right)+\boldsymbol{f}(\boldsymbol{x}_i^{(k)})\right\}
\right|\right|_2 \right]^2 \nonumber\\
= \min_{A \in \mathbb{R}^{m \times n} }\sum_{j=1}^N\left( d_{j(i)}^{(k)}   \left|\left| \Delta \boldsymbol{y}_{j(i)}^{(k)} - A \Delta \boldsymbol{x}_{j(i)}^{(k)}
\right|\right|_2 \right)^2 \label{eq::linAppSum}
\end{eqnarray}
for $ i=1,...,N $, where $d_{j(i)}^{(k)}\geq 0$, $ j=1,...,N $ are weights.
Here, $\Delta \boldsymbol{y}_{j(i)}^{(k)} = \boldsymbol{f}(\boldsymbol{x}_j^{(k)})-\boldsymbol{f}(\boldsymbol{x}_i^{(k)})\in \mathbb{R}^m$ and $\Delta \boldsymbol{x}_{j(i)}^{(k)} = \boldsymbol{x}_j^{(k)}-\boldsymbol{x}_i^{(k)}\in \mathbb{R}^n$.  (Note that $\Delta \boldsymbol{y}_{i(i)}^{(k)} = \boldsymbol{0}$, $\Delta \boldsymbol{x}_{i(i)}^{(k)} = \boldsymbol{0}$.)
Also, let
\begin{eqnarray}
	\Delta Y_{(i)}^{(k)}&=&\left[\Delta \boldsymbol{y}_{1(i)}^{(k)},\Delta \boldsymbol{y}_{2(i)}^{(k)}, \dots \Delta \boldsymbol{y}_{N(i)}^{(k)} \right] \in \mathbb{R}^{m \times N}\\
	\Delta X_{(i)}^{(k)}&=&\left[\Delta \boldsymbol{x}_{1(i)}^{(k)},\Delta \boldsymbol{x}_{2(i)}^{(k)}, \dots \Delta \boldsymbol{x}_{N(i)}^{(k)} \right] \in \mathbb{R}^{n \times N}.
\end{eqnarray}
Note that $\boldsymbol{f}(\boldsymbol{x}_{i}^{(k)})$ are always computed at the previous iteration (e.g., as equation~\eqref{eq:initMatrixY} when $k=0$ and in Step 2-3 when $k>0$). Hence, no new evaluation of $\boldsymbol{f}$ is required at this step.

The key idea here is that we weight the information of the function evaluation near $\boldsymbol{x}_{i}^{(k)}$ more than the function evaluation further away. That is to say, $d_{j(i)}^{(k)}>d_{j'(i)}^{(k)}$ if $||\boldsymbol{x}_{j}^{(k)}-\boldsymbol{x}_{i}^{(k)}||<||\boldsymbol{x}_{j'}^{(k)}-\boldsymbol{x}_{i}^{(k)}||$.  The importance of this idea can be seen in the numerical experiment presented in Appendix~\ref{app::NumExpWeight}.

Noting that
\begin{eqnarray}
	\sum_{j=1}^N\left(d_{j(i)}^{(k)}\left|\left|A\Delta\boldsymbol{x}_{j(i)}^{(k)}-\Delta\boldsymbol{y}_{j(i)}^{(k)} \right|\right|_2 \right)^2&=&\left|\left|\left(A\Delta X_{(i)}^{(k)} -\Delta Y_{(i)}^{(k)} \right)D^{(k)}_{(i)} \right|\right|_\textrm{F}^{\,2}\\
	&=&\left|\left|D^{(k)}_{(i)}\left( {\Delta X_{(i)}^{(k)} }^\textrm{T} A^{\textrm{T}}- { \Delta Y_{(i)}^{(k)} }^\textrm{T} \right) \right|\right|_\textrm{F}^{\,2}\,,
\end{eqnarray}

 we can rewrite Equation~\eqref{eq::linAppSum} as
\begin{eqnarray}
\min_{A^{(k)}_{(i)}\in \mathbb{R}^{m\times n}} \left|\left|D^{(k)}_{(i)}\left( { \Delta X^{(k)}_{(i)} }^\textrm{T} A^{(k)\textrm{T}}_{(i)} - { \Delta Y^{(k)}_{(i)} }^\textrm{T} \right) \right|\right|^2_\textrm{F}	\label{eq::linApprox}
\end{eqnarray}
where 
\begin{eqnarray}
	D_{(i)}^{(k)}&=&\textrm{diag}\left( d^{(k)}_{1(i)},d^{(k)}_{2(i)},..., d^{(k)}_{N(i)}\right),
\end{eqnarray}
where $d^{(k)}_{l(i)}\geq0$.
In this paper, we choose the weights as
\begin{eqnarray}d_{j(i)}^{(k)}=\left\{ 
\begin{array}{ll}\left(\frac{1}{
\sum_{l=1}^{n}(({x}_{lj}^{(k)}-{x}_{li}^{(k)})/(x_l^\textrm{U}-x_l^\textrm{L}))^2}\right)^\gamma & \textrm{if }j\neq i\\
0&\textrm{if }j = i
\end{array}
\right. ,\label{eq::weights}
\end{eqnarray}
where ${x}_{lj}^{(k)}, x_l^\textrm{U}, x_l^\textrm{L}$ are the $l$ th element of the vectors $\boldsymbol{x}_{j}^{(k)}, \boldsymbol{x}^\textrm{U}, \boldsymbol{x}^\textrm{L}$, respectively ($l=1,...,n$), and $\gamma \geq 0$ is a constant.
We use this weighting scheme so that the ``information'' from the
nonlinear function evaluation from the point closer to the point of
approximation is more influential when constructing the linear
approximation. The distance between $\boldsymbol{x}_{i}$ and
$\boldsymbol{x}_{j}$ are normalised by the size of the domain of initial guess (i.e., $\boldsymbol{x}^\textrm{U}$ and
$\boldsymbol{x}^\textrm{L}$).  The effect of the weight $d_{j(i)}^{(k)}$ and the parameter $\gamma$ and its necessity is analysed in Appendix~\ref{app::NumExpWeight}. 
The minimum norm solution of Equation~\eqref{eq::linApprox} is given by
\begin{eqnarray}
	A^{(k)}_{(i)} = \Delta Y^{(k)}_{(i)} D^{(k)}_{(i)} \left( \Delta X^{(k)}_{(i)} D^{(k)}_{(i)} \right)^\dagger .
\end{eqnarray}
where $^\dagger$ denotes the Moore-Penrose inverse.

If $ {\rm rank} \left(\Delta X^{(k)}_{(i)} D^{(k)}_{(i)}\right) = n$,
\begin{eqnarray}
  A^{(k)}_{(i)} = \Delta Y^{(k)}_{(i)} D^{(k)}_{(i)} \left( \Delta X^{(k)}_{(i)} D^{(k)}_{(i)} \right)^{\textrm{T}} \left\{ \left( \Delta X^{(k)}_{(i)} D^{(k)}_{(i)} \right) \left( \Delta X^{(k)}_{(i)} D^{(k)}_{(i)} \right)^{\textrm{T}} \right\}^{-1} .
\end{eqnarray}
Generically, $\textrm{rank} \Delta X^{(k)}_{(i)} D^{(k)}_{(i)} = n.$ \hspace{2mm} $ \textrm{rank} \Delta X^{(k)}_{(i)} < n $ can happen when $ x^{(k)}_{li} =  c \hspace{3mm} (i=1,2,\ldots, N)$, i.e. when 
$\boldsymbol{x}^{(k)}_{i}$ lie in the same hyperplane $x_l = c$. This happens when the $l$-th component of  $\boldsymbol{x}^{(k)}_i \hspace{3mm} (i=1,2,\ldots,N)$ are all equal.

\YA{The concise pseudo-code for the weighted linear approximation can be found in Algorithm~\ref{alg:linApprox}.}

\begin{algorithm}     
                   
\caption{$A^{(k)}_{(i)}$={\bf construct\_linearApproximation}$(i,X^{(k)},Y^{(k)},\boldsymbol{x}^\textrm{L},\boldsymbol{x}^\textrm{U},\gamma)$}         
\label{alg:linApprox}                          
\begin{algorithmic}                  
\STATE{$D^{(k)}_{(i)} \gets$ allocate memory for diagonal matrix of $N \times N$}
\STATE{$\Delta X^{(0)}\gets$ allocate memory for $n \times N$ matrix}
\STATE{$\Delta Y^{(0)}\gets$ allocate memory for $m \times N$ matrix}

\FORALL{$j= 1,...,N$}
\STATE{$d_{jj(i)}^{(k)}\gets \left\{ 
\begin{array}{ll}\left(\frac{1}{
\sum_{l=1}^{n}(({x}_{lj}^{(k)}-{x}_{li}^{(k)})/(x_l^\textrm{U}-x_l^\textrm{L}))^2}\right)^\gamma & \textrm{if }j\neq i\\
0&\textrm{if }j = i
\end{array}
\right. $}
\STATE{$\Delta\boldsymbol{x}^{(k)}_j\gets \boldsymbol{x}^{(k)}_j-\boldsymbol{x}^{(k)}_i$}
\STATE{$\Delta\boldsymbol{y}^{(k)}_j\gets \boldsymbol{y}^{(k)}_j-\boldsymbol{y}^{(k)}_i$}
\ENDFOR
\STATE{$A^{(k)}_{(i)} \gets \Delta Y^{(k)}_{(i)} D^{(k)}_{(i)} \left( \Delta X^{(k)}_{(i)} D^{(k)}_{(i)} \right)^\dagger $} \COMMENT{$\dagger$: Moore-Penrose inverse.}
\end{algorithmic}
\end{algorithm}

\textbf{2-2) Solve for $\boldsymbol{x}$ that
minimises
$||\boldsymbol{y}^*- (A_{(i)}^{(k)} (\boldsymbol{x}-\boldsymbol{x}_{i}^{(k)})+\boldsymbol{f}(\boldsymbol{x}_{i}^{(k)}))||^{\,2}_2$}\\
We now compute the next iterate $X^{(k+1)}$ using the matrices
$\{A^{(k)}_{(i)}\}_{i=1}^{N}$ similarly to the
Gauss-Newton method with Tikhonov regularisation (e.g., \cite{hansen2005rank, Bjorck1996}), i.e.,

\begin{equation}\boldsymbol{x}_{i}^{(k+1)}=\boldsymbol{x}_{i}^{(k)}+\left(A_{(i)}^{(k)\textrm{T}}A_{(i)}^{(k)}+\lambda_i^{(k)} I\right)^{-1}A_{(i)}^{(k)\textrm{T}}(\boldsymbol{y}^*-\boldsymbol{y}_{i}^{(k)})\qquad  \textrm{for }i=1,...,N, \label{eq::regularise}
\end{equation}
where $\boldsymbol{y}^*$ is the set of observations one wishes to fit the nonlinear function $\boldsymbol{f}$ to (cf. \eqref{eq:mathProblem}), and $\boldsymbol{y}_i^{(k)}\equiv \boldsymbol{f}(\boldsymbol{x}_i^{(k)})$.  For CGN, we require $\lambda_i^{(k)}> 0$. The necessity of the regularisation can be seen in the numerical experiment presented in Appendix~\ref{app::NumExpReg}.\\

\textbf{2-3) Update matrices $X$ and $Y$ and vectors $\boldsymbol{r}$ and
$\boldsymbol{\lambda}$}\\
Evaluate the nonlinear function $\boldsymbol{f}$ for each
$\boldsymbol{x}_{i}^{(k+1)}$ as $\boldmath{y}_i^{(k+1)} = \boldmath{f}(\boldsymbol{x}_i^{(k+1)} )\hspace{2mm} (i=1,2,\ldots,n)$, and store as $Y^{(k+1)}= [ \boldmath{y}_1^{(k+1)}, \boldmath{y}_2^{(k+1)},\ldots,\boldsymbol{y}_N^{(k+1)} ].$
Compute the sum of squared residuals vector as $\boldsymbol{r}^{(k+1)}= [ r_1^{(k+1)}, r_2^{(k+1)},\ldots, r_N^{(k+1)} ]^\textrm{T}$ 
where $ r_i^{(k+1)} = \| \boldmath{y}_i^{(k+1)} - \boldmath{y}^{\ast} \|_2^{\,2} \hspace{2mm} (i=1,2,\ldots,N)$.
Note that this process can be implemented embarrassingly parallelly.\\

If the residual $r_i$ increases, we replace $\boldsymbol{x}_{i}^{(k+1)}$ by $\boldsymbol{x}_{i}^{(k)}$ and increase the regularisation parameter  i.e.,

if $r_{i}^{(k)}<r_{i}^{(k+1)}$ or
$\boldsymbol{f}(\boldsymbol{x}_{i}^{(k+1)})$ can not be
evaluated, then let

\begin{eqnarray}\boldsymbol{x}_{i}^{(k+1)}=\boldsymbol{x}_{i}^{(k)}\\
\boldsymbol{y}_{i}^{(k+1)}=\boldsymbol{y}_{i}^{(k)}\\
\lambda^{(k+1)}_i=10 \,\lambda^{(k)}_i\end{eqnarray}

else decrease the regularisation parameter, i.e.,

\begin{eqnarray}\lambda^{(k+1)}_i= \frac{1}{10}\lambda^{(k)}_i \,,\end{eqnarray}

\noindent for each $i=1,..., N$.  

\YA{There are various ways to update the regularisation parameter $\lambda$. In this paper we followed the strategy used in the Matlab's implementation of the Levenberg-Marquardt method.}

\YA{In addition, in this step, we impose the stopping criteria for each point in the iteration. 
As can be seen in Equation $\eqref{eq::regularise}$, $\boldsymbol{x}_{i}^{(k+1)}\approx\boldsymbol{x}_{i}^{(k)}$ for large $\lambda_i$, so that we can expect very small update in $\boldsymbol{x}_{i}^{(k+1)}$. Hence, in order to mimic the minimum step size stopping criteria, we stop the update 
for $i$ where $\lambda_i>\lambda_\textrm{max}$.}

\YA{The concise pseudo-code for updating matrices $X$ and $Y$ and vectors $\boldsymbol{r}$ and
$\boldsymbol{\lambda}$ is given in Algorithm~\ref{alg:updateCluster}.}

The influence of the choice of the initial value $\lambda_{\textrm{init}}$ of the regularization parameter in equation (\ref{lambda_init}) is studied in Appendix D.

\setlength{\leftskip}{0cm}
}

\begin{algorithm}     
                   
\caption{\\$(X^{(k+1)},Y^{(k+1)},\boldsymbol{r}^{(k+1)},\boldsymbol{\lambda}^{(k+1)})$={\bf update\_cluster}$(X^{(k)},X^{(k+1)},Y^{(k)},\boldsymbol{r}^{(k)},\boldsymbol{\lambda}^{(k)},\lambda_\textrm{\tiny max})$}         
\label{alg:updateCluster}                          
\begin{algorithmic}
\STATE{$Y^{(k+1)}\gets$ allocate memory for $m \times N$ matrix}
\STATE{$\boldsymbol{r}^{(k+1)}\gets$ allocate memory for length $N$ vector}
\STATE{$\boldsymbol{\lambda}^{(k+1)}\gets$ allocate memory for length $N$ vector}
                  
\FORALL{$i= 1,...,N$} 
\IF{$\lambda_i^{(k)}\leq\lambda_\textrm{max}$ and $\boldsymbol{f}(\boldsymbol{x}_i^{(k+1)})$ can be evaluated}
\STATE{$\boldsymbol{y}_i^{(k+1)}\gets \boldsymbol{f}(\boldsymbol{x}_i^{(k+1)})$ }
\STATE{$r_i^{(k+1)} \gets || \boldsymbol{y}_i^{(k+1)}-\boldsymbol{y}^*||_2^{\,2}$}

\IF{$r_i^{(k+1)}>r_i^{(k)}$}
\STATE{$\boldsymbol{x}_i^{(k+1)}\gets \boldsymbol{x}_i^{(k)}$ }
\STATE{$\boldsymbol{y}_i^{(k+1)}\gets \boldsymbol{y}_i^{(k)}$ }
\STATE{$r_i^{(k+1)} \gets r_i^{(k)}$}
\STATE{$\lambda_i^{(k+1)}\gets 10\,\lambda_i^{(k)}$}
\ELSE
\STATE{$\lambda_i^{(k+1)}\gets \frac{1}{10}\lambda_i^{(k)}$}
\ENDIF

\ELSE
\STATE{$\boldsymbol{x}_i^{(k+1)}\gets \boldsymbol{x}_i^{(k)}$ }
\STATE{$\boldsymbol{y}_i^{(k+1)}\gets \boldsymbol{y}_i^{(k)}$ }
\STATE{$r_i^{(k+1)} \gets r_i^{(k)}$}
\STATE{$\lambda_i^{(k+1)}\gets 10\,\lambda_i^{(k)}$}

\ENDIF
\ENDFOR

\end{algorithmic}
\end{algorithm}

\YA{We present a concise description of the CGN as a pseudo-code in Algorithm~\ref{alg:CGN}.  }

\begin{algorithm}
                   
\caption{\\$X^{(k_\textrm{max})}$ ={\bf Cluster\_Gauss-Newton\_method}$(\boldsymbol{x}^\textrm{L},\boldsymbol{x}^\textrm{U}, N, \boldsymbol{f},\boldsymbol{y}^*, \lambda_\textrm{init}, \lambda_\textrm{max}, \gamma, k_\textrm{max})$}         
\label{alg:CGN}                          
\begin{algorithmic}
\STATE{$\boldsymbol{r}^{(0)}\gets$ allocate memory for length $N$ vector}
\FORALL{$i \gets 1,N$}
\STATE{$\lambda_i^{(0)}=\lambda_\textrm{init}$}
\ENDFOR
\\
\STATE{$(X^{(0)},Y^{(0)},\boldsymbol{r}^{(0)})$={\bf create\_initialCluster}$(\boldsymbol{x}^\textrm{L},\boldsymbol{x}^\textrm{U}, N, \boldsymbol{f},\boldsymbol{y}^*)$}
\\
\FOR{$k \gets 0,(k_\textrm{max}-1)$}

\FORALL{$i \gets 1,N$}

\IF{$\lambda_i^{(k)}\leq\lambda_\textrm{max}$}
\STATE{$A^{(k)}_{(i)}$={\bf construct\_linearApproximation}$(i,X^{(k)},Y^{(k)},\boldsymbol{x}^\textrm{L},\boldsymbol{x}^\textrm{U},\gamma)$}
\STATE{$\boldsymbol{x}_{i}^{(k+1)}=\boldsymbol{x}_{i}^{(k)}+\left(A_{(i)}^{(k)\textrm{T}}A_{(i)}^{(k)}+\lambda_i^{(k)} I\right)^{-1}A_{(i)}^{(k)\textrm{T}}(\boldsymbol{y}^*-\boldsymbol{y}_{i}^{(k)})$}
\ENDIF

\ENDFOR
\STATE{$(X^{(k+1)},Y^{(k+1)},\boldsymbol{r}^{(k+1)},\boldsymbol{\lambda}^{(k+1)})$\\\qquad \qquad ={\bf update\_cluster}$(X^{(k)},X^{(k+1)},Y^{(k)},\boldsymbol{r}^{(k)},\boldsymbol{\lambda}^{(k)},\lambda_\textrm{ max})$}

\ENDFOR

\end{algorithmic}
\end{algorithm}

\section{Numerical Experiments}\label{header-c1748}
In this section we illustrate the advantages of the proposed CGN algorithm by numerical experiments on three PBPK models.
\subsection{Numerical Experiment setup}\label{header-c1749}
In this subsection we specify the details of the numerical experiment set-up.
\subsubsection{Mathematical models}\label{header-c1750}
 For the numerical experiments, we used the following three published mathematical models of drug concentration in the blood of a human body (PBPK model). The time course drug concentration was simulated using the model, and random noise was added to mimic the observation uncertainties. The random noise was generated by a normal distribution with a standard deviation of 10\% of the simulated concentration value. The simulated drug concentration was used as the \YA{test dataset}, and multiple possible parameters were estimated from \YA{the test dataset} using CGN and conventional methods. The implementations of all the examples are available as the supplementary material. 

{\bf Example 1: Multi-dose problem}\\
For this example, we consider the case where three different amounts of a drug is given (low-dose, medium-dose, and high-dose) orally as pills to a patient.  This can be modelled using the same mathematical model as in the Motivating Example (cf. Appendix~\ref{sec::motivatingExample} and Appendix~\ref{app::PBPKmodel}).  The initial value problem of this mathematical model can be written as:
\begin{eqnarray}
\dee{\boldsymbol{u}}{t}&=&\boldsymbol{g}(\boldsymbol{u},t;\boldsymbol{x})\,,\label{eq::Example1_ODE}\\
u_{i}(t=0)&=&0\qquad \textrm{for } i=1,...,17,19,20\,,\\
u_{18}(t=0)&=&\left\{\begin{array}{ll}30000&\textrm{for low-dose,}\\ 100000&\textrm{for medium-dose,} \\ 300000&\textrm{for high-dose.}\end{array}\right. \label{eq::Example1_IC}
\end{eqnarray}
In this mathematical model, $u_1(t)$ represents the drug concentration in blood at time $t$.  For convenience, we denote the drug concentration at time $t$ for the low-dose, medium-dose, and high-dose as $u^\textrm{l}_1(t;\boldsymbol{x})$, $u^\textrm{m}_1(t;\boldsymbol{x})$, and $u^\textrm{h}_1(t;\boldsymbol{x})$, respectively.  Note that the right hand side of the system of ODEs~\eqref{eq::Example1_ODE} depends on the parameter vector $\boldsymbol{x}$. Hence, the solution of the system of ODEs $\boldsymbol{u}$ depends not only on $t$ but also on $\boldsymbol{x}$.  For this example, we consider the case where the blood sample is taken at $t=2,3,4,6,8,12,24,36,48,72$, so that the nonlinear function $\boldsymbol{f}$ in \eqref{eq:mathProblem} can be written as
\begin{eqnarray}
\boldsymbol{f}(\boldsymbol{x})&=&[u^\textrm{l}_1(2;\boldsymbol{x}), u^\textrm{l}_1(3;\boldsymbol{x}), u^\textrm{l}_1(4;\boldsymbol{x}), u^\textrm{l}_1(6;\boldsymbol{x}), u^\textrm{l}_1(8;\boldsymbol{x}), u^\textrm{l}_1(12;\boldsymbol{x}), u^\textrm{l}_1(24;\boldsymbol{x}), u^\textrm{l}_1(36;\boldsymbol{x}), \nonumber\\
&&u^\textrm{l}_1(48;\boldsymbol{x}), u^\textrm{l}_1(72),u^\textrm{m}_1(2;\boldsymbol{x}), u^\textrm{m}_1(3;\boldsymbol{x}), u^\textrm{m}_1(4;\boldsymbol{x}), u^\textrm{m}_1(6;\boldsymbol{x}), u^\textrm{m}_1(8;\boldsymbol{x}), u^\textrm{m}_1(12;\boldsymbol{x}), \nonumber\\
&&u^\textrm{m}_1(24;\boldsymbol{x}), u^\textrm{m}_1(36;\boldsymbol{x}), u^\textrm{m}_1(48;\boldsymbol{x}), u^\textrm{m}_1(72),u^\textrm{h}_1(2;\boldsymbol{x}), u^\textrm{h}_1(3;\boldsymbol{x}), u^\textrm{h}_1(4;\boldsymbol{x}), u^\textrm{h}_1(6;\boldsymbol{x}), \nonumber\\
&&u^\textrm{h}_1(8;\boldsymbol{x}), u^\textrm{h}_1(12;\boldsymbol{x}), u^\textrm{h}_1(24;\boldsymbol{x}), u^\textrm{h}_1(36;\boldsymbol{x}), u^\textrm{h}_1(48;\boldsymbol{x}), u^\textrm{h}_1(72)]^\textrm{T}.\nonumber
\end{eqnarray}
This example is based on \cite{Fukuchi2017}.

{\bf Example 2: Intravenous injection and oral administration problem}\\
For this example, we consider the case where the drug is administered via injection into the vein (i.v.) and the case where the drug is given orally as a pill (p.o.) to a patient.  This was modelled using a mathematical model similar to the motivating example in Appendix~\ref{sec::motivatingExample}.  The initial value problem of this mathematical model for i.v. administration can be written as:
\begin{eqnarray}
\dee{\boldsymbol{u}}{t}&=&\boldsymbol{\tilde{g}}(\boldsymbol{u},t;\boldsymbol{x})\,,\label{eq::Example2_ODE}\\
u_{1}(t=0)&=&\frac{30.488}{0.074+10^{x_4}}\,,\\
u_{i}(t=0)&=&0\qquad \textrm{for } i=2,3,...,20\,,
\end{eqnarray}
where $x_4$ is one of the parameters in the parameter vector $\boldsymbol{x}$.  The initial value problem of this mathematical model for p.o. administration can be written as:
\begin{eqnarray}
\dee{\boldsymbol{u}}{t}&=&\boldsymbol{\tilde{g}}(\boldsymbol{u},t;\boldsymbol{x})\,,\label{eq::Example22_ODE}\\
u_{i}(t=0)&=&0\qquad \textrm{for } i=1,...,20\,,\\
u_{18}\left(t=\frac{e^{x_{11}}}{2(1+e^{x_{11}})}\right)&=&30.488\,,\\
u_{i}\left(t=\frac{e^{x_{11}}}{2(1+e^{x_{11}})}\right)&=&0\qquad \textrm{for } i=1,...,17,19,20\,,
\end{eqnarray}
where $x_{11}$ is one of the parameters in the parameter vector $\boldsymbol{x}$. $\frac{e^{x_{11}}}{2(1+e^{x_{11}})}$ is the delay in the absorption of the drug (for example the time it takes from intake in the mouth to pill being dissolved in the stomach).  Similarly to Example~1, the observable quantity, the drug concentration in the blood plasma is represented as $u_1(t;\boldsymbol{x})$. For convenience, we denote the drug concentration at time $t$ when the drug is administrated by i.v. and p.o. administrations as $u^\textrm{i}_1(t;\boldsymbol{x})$ and $u^\textrm{p}_1(t;\boldsymbol{x})$, respectively.  For this example, we consider the case where the blood sample is taken at $t=0.0833, 0.1667, 0.25, 0.5, 0.75, 1, 1.5, 2, 3, 4, 6, 8$ when the drug is given as intravenous injection and the blood sample is taken at $t=0.5, 1, 1.5, 2, 3, 4, 6, 8, 12, 14$ when the drug is given as a pill orally.  The nonlinear function $\boldsymbol{f}$ can be written as
\begin{eqnarray}
\boldsymbol{f}(\boldsymbol{x})&=&[u_1^i(0.0833;\boldsymbol{x}),u_1^i(0.1667;\boldsymbol{x}),u_1^i(0.25;\boldsymbol{x}),u_1^i(0.5;\boldsymbol{x}),u_1^i(0.75;\boldsymbol{x}),u_1^i(1;\boldsymbol{x}),\nonumber\\
&&u_1^i(1.5;\boldsymbol{x}),u_1^i(2;\boldsymbol{x}),u_1^i(3;\boldsymbol{x}),u_1^i(4;\boldsymbol{x}),u_1^i(6;\boldsymbol{x}),u_1^i(8;\boldsymbol{x}),u_1^p(0.5;\boldsymbol{x}),u_1^p(1;\boldsymbol{x}),\nonumber\\
&&u_1^p(1.5;\boldsymbol{x}),u_1^p(2;\boldsymbol{x}),u_1^p(3;\boldsymbol{x}),u_1^p(4;\boldsymbol{x}),u_1^p(6;\boldsymbol{x}),u_1^p(8;\boldsymbol{x}),u_1^p(12;\boldsymbol{x}),u_1^p(14;\boldsymbol{x})]^\textrm{T}. \nonumber
\end{eqnarray}
This example is based on \cite{yao2018quantitative}. \YA{For this example, we used the upper and lower bounds of the initial cluster of CGN, which we know to include the known global minimiser (based on the parameter set we simulated the data from), as the bound constraints for the constrained optimisation algorithms. However, in practice, it is usually very hard to specify the upper and lower bounds of the parameters for the parameter estimation problem of our interest.}
 
{\bf Example 3: Drug-drug interaction problem}\\
For this example, we consider the case when a patient takes two different drugs.  As is often stated in the instruction for drugs, if two drugs are taken together, they can potentially interact inside the body and cause undesirable effects.  We model the cases where pitavastatin (for the ease of writing we shall refer to this drug as Drug A) is taken alone or with cyclosporin A (Drug B).  The concentration of Drug A is modelled using a mathematical model similar to the Motivating Example (Appendix~\ref{sec::motivatingExample}) and Drug B is modelled using a simplified version of the model. The interaction of Drugs A and B in the liver compartment, when they are administered simultaneously is also modelled.  

The initial value problem of the mathematical model for Drug A can be written as
\begin{eqnarray}
\dee{\boldsymbol{u}}{t}&=&\boldsymbol{g}(\boldsymbol{u},t;\boldsymbol{x})\\
u_{i}(t=0)&=&0\qquad \textrm{for } i=1,...,20\,,\\
u_{12}\left(t=\frac{e^{x_{18}}}{2(1+e^{x_{18}})}\right)&=&30.4971\,,\\
u_{i}\left(t=\frac{e^{x_{18}}}{2(1+e^{x_{18}})}\right)&=&0\qquad \textrm{for } i=1,...,11,13,...,20\,.
\end{eqnarray}

The initial value problem of the mathematical model for Drug A administered with Drug B can be written as:
\begin{eqnarray}
\dee{\boldsymbol{u}}{t}&=&\boldsymbol{h}(\boldsymbol{u},t;\boldsymbol{x})\\
u_{i}(t=0)&=&0\qquad \textrm{for } i=1,...,33\,,\\
u_{12}\left(t=\frac{e^{x_{18}}}{2(1+e^{x_{18}})}\right)&=&30.4971\,,\\
u_{33}\left(t=\frac{e^{x_{19}}}{2(1+e^{x_{19}})}\right)&=&2000\,,\\
u_{i}\left(t=\frac{e^{x_{18}}}{2(1+e^{x_{18}})}\right)&=&0\qquad \textrm{for } i=1,...,11,13,..,33\,,\\
u_{i}\left(t=\frac{e^{x_{19}}}{2(1+e^{x_{19}})}\right)&=&0\qquad \textrm{for } i=1,...,32\,.
\end{eqnarray}
For convenience, we denote the concentration of Drug A at time $t$ when only Drug A is administered as $u^\textrm{A}_1(t;\boldsymbol{x})$.  When both Drugs A and B are administered, we denote the drug concentration of Drug A  as $u^\textrm{AB}_1(t;\boldsymbol{x})$ and the concentration of Drug B as $u^\textrm{AB}_{21}(t;\boldsymbol{x})$.
For this example, we consider the case where the blood sample is taken at $t=0.5,1,1.5,2,3,5,8,12$, and for the case where both Drugs A and B are administered, we measure the drug concentration of both drugs.  The nonlinear function $\boldsymbol{f}$ can be written as
\begin{eqnarray}
\boldsymbol{f}(\boldsymbol{x})&=&[u^\textrm{A}_1(0.5;\boldsymbol{x}), u^\textrm{A}_1(1;\boldsymbol{x}), u^\textrm{A}_1(1.5;\boldsymbol{x}), u^\textrm{A}_1(2;\boldsymbol{x}), u^\textrm{A}_1(3;\boldsymbol{x}), u^\textrm{A}_1(5;\boldsymbol{x}), u^\textrm{A}_1(8;\boldsymbol{x}),\nonumber\\
&& u^\textrm{A}_1(12;\boldsymbol{x}),u^\textrm{AB}_1(0.5;\boldsymbol{x}), u^\textrm{AB}_1(1;\boldsymbol{x}), u^\textrm{AB}_1(1.5;\boldsymbol{x}), u^\textrm{AB}_1(2;\boldsymbol{x}), u^\textrm{AB}_1(3;\boldsymbol{x}), u^\textrm{AB}_1(5;\boldsymbol{x}),\nonumber\\
&& u^\textrm{AB}_1(8;\boldsymbol{x}), u^\textrm{AB}_1(12;\boldsymbol{x}), u^\textrm{AB}_{21}(0.5;\boldsymbol{x}), u^\textrm{AB}_{21}(1;\boldsymbol{x}), u^\textrm{AB}_{21}(1.5;\boldsymbol{x}), u^\textrm{AB}_{21}(2;\boldsymbol{x}), u^\textrm{AB}_{21}(3;\boldsymbol{x}), \nonumber\\
&&u^\textrm{AB}_{21}(5;\boldsymbol{x}), u^\textrm{AB}_{21}(8;\boldsymbol{x}), u^\textrm{AB}_{21}(12;\boldsymbol{x})]^\textrm{T}\nonumber
\end{eqnarray}
This example is based on \cite{yoshikado2016quantitative}.

The mathematical description of each PBPK model used for the numerical experiments is summarised in Table~1.

\begin{table}
\begin{longtable}[]{@{}llll@{}}
\toprule
& \multirow{ 2}{*}{Model structure} & Parameters & Simulated
\tabularnewline
&  & to be Estimated & Observations\tabularnewline
\midrule
\endhead
\multirow{ 2}{*}{Model 1} & Three systems of nonlinear ODEs& \multirow{ 2}{*}{11 parameters} & \multirow{ 2}{*}{30 observations} \tabularnewline
&with 20 variables & & \tabularnewline
\tabularnewline
\multirow{ 2}{*}{Model 2} & Two systems of linear ODEs& \multirow{ 2}{*}{11 parameters} & \multirow{ 2}{*}{22 observations}\tabularnewline
&with 20 variables & &\tabularnewline
\tabularnewline
\multirow{ 2}{*}{Model 3} & Two systems of nonlinear ODEs& \multirow{ 2}{*}{19 parameters} & \multirow{ 2}{*}{24 observations}\tabularnewline
&with 20 and 33 variables & & \tabularnewline
%
%
\bottomrule
\end{longtable}
\caption{Summary of the mathematical description of the PBPK models used for the numerical
experiments.}
\end{table}

\subsubsection{Computation environment}\label{header-c1777}

All computational experiments were performed using Matlab \YA{2019a} on 3.1
GHz Intel Core i5 processors with MacOS version 10.\YA{13.6}. \YA{When using the algorithms implemented in Python we used Python version 3.7.}  All results of the numerical experiments were summarised and visualised using ggplot2 version 2.2.1 \cite{ggplot} in R version 3.3.2.

\subsubsection{The initial set of vectors
$\{\boldsymbol{x}_{i}^{(0)}\}_{i=1}^{N}$} \label{sec::initialVector} \YA{We generated the initial set of vectors $\{\boldsymbol{x}_{i}^{(0)}\}_{i=1}^{N}$ randomly based on Algorithm~\ref{alg:initialCluster}. The range of the initial cluster was set by the pharmacologically likely parameter range based on a priori knowledge defined by domain specialists} (e.g., from the values obtained from the animal experiments or lab experiments).  

\YA{We stored and used this initial set of vectors as the initial cluster for the CGN and also as initial iterates for the other nonlinear least squares solvers and optimisation algorithms which were compared.}

\YA{Notice that this range of the initial cluster does not necessary contain the global minimiser as many pharmacokinetics parameter can be few order of magnitudes different between species.  For the sake of comparison with the constrained optimisation algorithms, we constructed Example 2 so that the initial range contains the set of parameters that we used to create the dataset.  This ensures that a global minimiser is in the domain where the initial cluster is made.  When comparing with the constrained optimisation algorithms, we used this domain as the bound constraint for the parameter search for these algorithms.}

\subsubsection{ODE solver}\label{header-c1784} For Examples 1 and 3 the nonlinear function evaluations require the numerical solution of stiff systems of ODEs. We used the ODE15s solver \cite{Shampine1997} with the default settings to solve these ODEs \YA{(relative tolerance $10^{-3}$, absolute tolerance $10^{-6}$)}. We observed that for some set of parameters, the ODE solver can get stuck in an infinite loop. Here, we set the timeout, where if the ODE evaluation takes longer than 5 seconds, it terminates the solution process and returns a not-a-number vector.  

\subsubsection{Setting for the Cluster Gauss-Newton (CGN) method:}\label{header-c1792} We used the following parameters unless stated otherwise:d
\begin{eqnarray}
N&=&250\\
\lambda_\textrm{init}&=&0.01\\
\lambda_\textrm{max}&=&10^{10}\\
k_\textrm{max}&=&100.
\end{eqnarray}
Here $k_\textrm{max}$ is the maximum number of iterations of the method. 
We used the initial set of vectors $\{\boldsymbol{x}_{i}^{(0)}\}_{i=1}^{N}$ described in Section~\ref{sec::initialVector} as the initial cluster. \YA{Note that we chose $\lambda_\textrm{init}$ to be consistent with the default setting of the LM implementation (\emph{lsqnonlin} in Matlab.}
 
\subsection{Algorithms Compared\label{header-c1780}}
\YA{Here, we list the conventional and recently developed algorithms that we compared our proposed algorithm with.  We used the default setting of the algorithms unless stated otherwise.}

\subsubsection{Nonlinear least squares solvers}
\YA{We compare the proposed algorithm against gradient-based and gradient-free nonlinear least squares algorithms.
We used these algorithms to find multiple approximate minimisers of the nonlinear least square problem of our interest by repeatedly applying these algorithms with various initial iterates. Namely; we used each parameter vector in $\{\boldsymbol{x}_{i}^{(0)}\}_{i=1}^{N}$, which were randomly generated in the CGN method as the initial iterate, and excecuted each algorithm repeatedly to obtain $N$ sets of estimated parameter vectors.}

{\bf Levenberg-Marquardt (LM) method:}\label{header-c1787} The conventional gradient-based nonlinear least squares solver, Levenberg-Marquardt method \cite{Bjorck1996} implemented in the \emph{lsqcurvefit} function in Matlab. 

We use LM with the default setting as well as setting `FiniteDifferenceStepSize' to be the square root of the default accuracy of the ODE solve (e.g., \\
$\textrm{FiniteDifferenceStepSize}=\sqrt{\textrm{AbsTol}}=\sqrt{10^{-6}}$).  Note that the\\
 FiniteDifferenceStepSize of the default setting is $10^{-6}$.

{\bf Trust-Region (TR) method:} The conventional gradient-based nonlinear least squares solver, Trust-Region method \cite{conn2000trust} implemented in the \emph{lsqcurvefit} function in Matlab.

{\bf DFO-LS method:}
\YA{A gradient-free nonlinear least squares solver DFO-LS \\
method\cite{cartis2018improving}.
We obtained the Python code of DFO-LS Version 1.0.2 from \url{http://people.maths.ox.ac.uk/robertsl/dfols/} (last accessed on June 6th 2019).  To make the numerical method for solving the model ODEs exactly the same, we called the Matlab implementation of the nonlinear function through MATLAB Engine API for Python.}


{\bf libensemble with POUNDERS:}
\YA{We used the libensemble algorithm \cite{libEnsemble} that was available at the developer branch of \url{https://github.com/Libensemble/libensemble} (last accessed April 25th 2019) together with petsc version 3.10.4 available as \url{https://www.mcs.anl.gov/petsc/index.html}.  We  used POUNDERS as the nonlinear least squares solver. As libensemble was implemented in Python, we used the MATLAB Engine API to call the nonlinear function. POUNDERS algorithm utilises bound constraints for the search domain, and we tested this algorithm only for Example~2.}
%
%
%
%
\subsubsection{Optimisation algorithms}
\YA{We compared the proposed method against the optimisation algorithms that are designed to minimise a function which takes a vector quantity as an input and scalar quantity as output.  We applied these algorithms to our nonlinear least squares problem by minimising SSR defined in Equation~\eqref{eq::SSR}}

{\bf Quasi-Newton method:}
\YA{We used the Quasi-Newton method implemented in Matlab's Optimisation toolbox as the \emph{fminunc} function.  This implementation uses the BFGS formula to update the approximate Hessian.}

The following algorithms require bound constraints, and they were tested only for Example~2.

{\bf Implicit Filter method:}
\YA{This is another gradient-free optimisation algorithm.  The code was downloaded from \url{https://archive.siam.org/books/se23/} (last accessed June 7th 2019).  We used the `least-squares' option with a budget of 200.}

{\bf Surrogate optimisation method:}
\YA{We used the Surrogate optimisation algorithm implemented in Matlab's Global Optimisation Toolbox as the \emph{surrogateopt} function. We repeated this algorithm 250 times with distinct random-seed to obtain 250 global minimisers of the nonlinear-least squares problem of our interest.}

{\bf Particle Swarm}
\YA{We used the Particle Swarm optimisation algorithm implemented in Matlab's Global Optimisation Toolbox as the \emph{particleswarm} function. We repeated this algorithm 250 times with distinct random seeds to obtain 250 global minimisers of the nonlinear least squares problem of our interest.}

{\bf Genetic Algorithm:}
\YA{We used the Genetic algorithm implemented in Matlab's Global Optimisation Toolbox as the \emph{ga} function. As one run of this algorithm required over 100,000 function evaluations for our examples, we did not repeat this algorithm to obtain multiple global minimisers.}

{\bf DIRECT:}
\YA{A sampling algorithm DIRECT \cite{jones1993lipschitzian}. The Matlab implementation was obtained from \url{https://searchcode.com/codesearch/view/12449743/} (last accessed July 30th 2019). As one run of this algorithm required over 100,000 function evaluations for our examples, we did not repeat this algorithm to obtain multiple global minimisers.}

\subsection{Results}
 \YA{We compared the proposed CGN method with various existing algorithms.  We compared the `speed' of the algorithm by the number of function evaluations required and we compared the `quality' of the minimisers found by the SSR (smaller the SSR the better minimiser). As we can see in Table~\ref{tab::compTimeCGN}, the dominant part of the computation time was spent by the nonlinear function evaluations in CGN. Hence, we claim that the number of function evaluations is a fair way to compare the computation cost.}

\begin{table}[ht]
\centering
\begin{tabular}{rrrr}
  \hline
 &Example 1&Example 2&Example 3 \\ 
  \hline
 Nonlinear function evaluation&1,890.598 & 852.116 & 1,779.362 \\ 
 Computation of linear approximation & 3.610 & 6.929 & 11.021 \\ 
 Output results as csv files & 3.647 & 3.663 &4.129 \\ 
 \hline
 Total computation & 1,898.393 & 863.172 & 1,795.192 \\
 \hline
\end{tabular}
\caption{Ingredients of the computation time for CGN (seconds).} 
\label{tab::compTimeCGN}
\end{table}

\begin{figure}

\begin{subfigure}[b]{1\columnwidth}
\includegraphics[width=110mm]{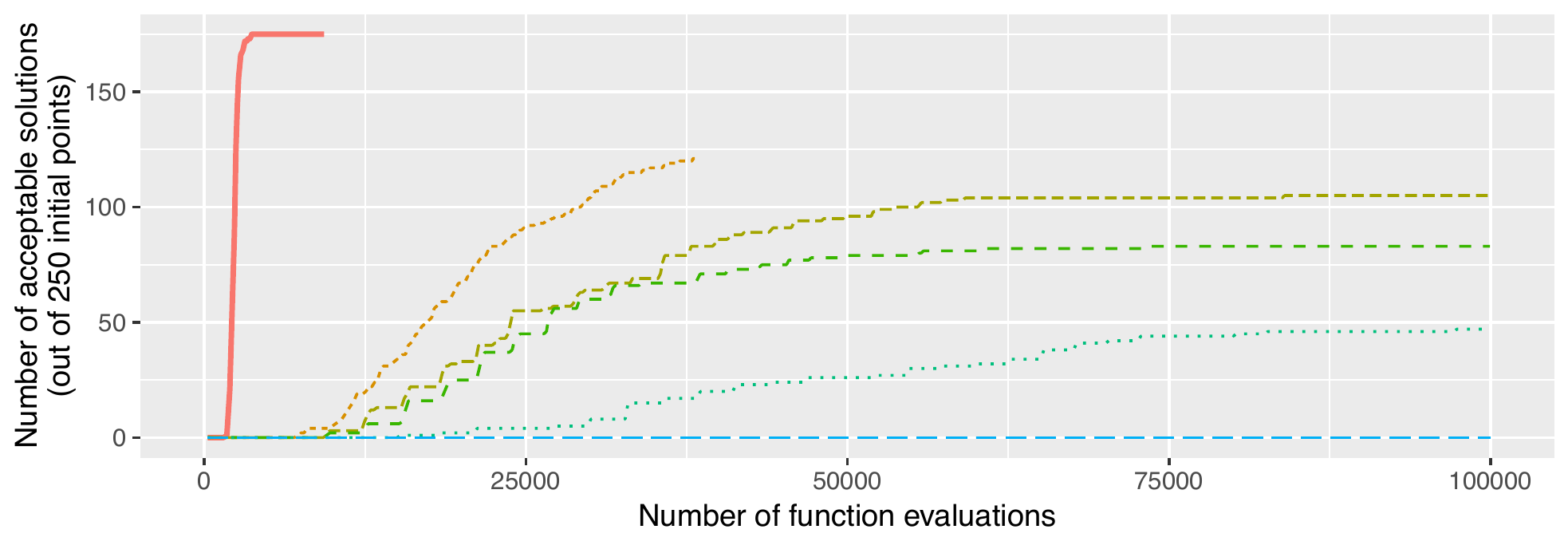}\label{fig::Example1_speed_plot}
\caption{Example 1}
\end{subfigure}
\begin{subfigure}[b]{1\columnwidth}
\includegraphics[width=110mm]{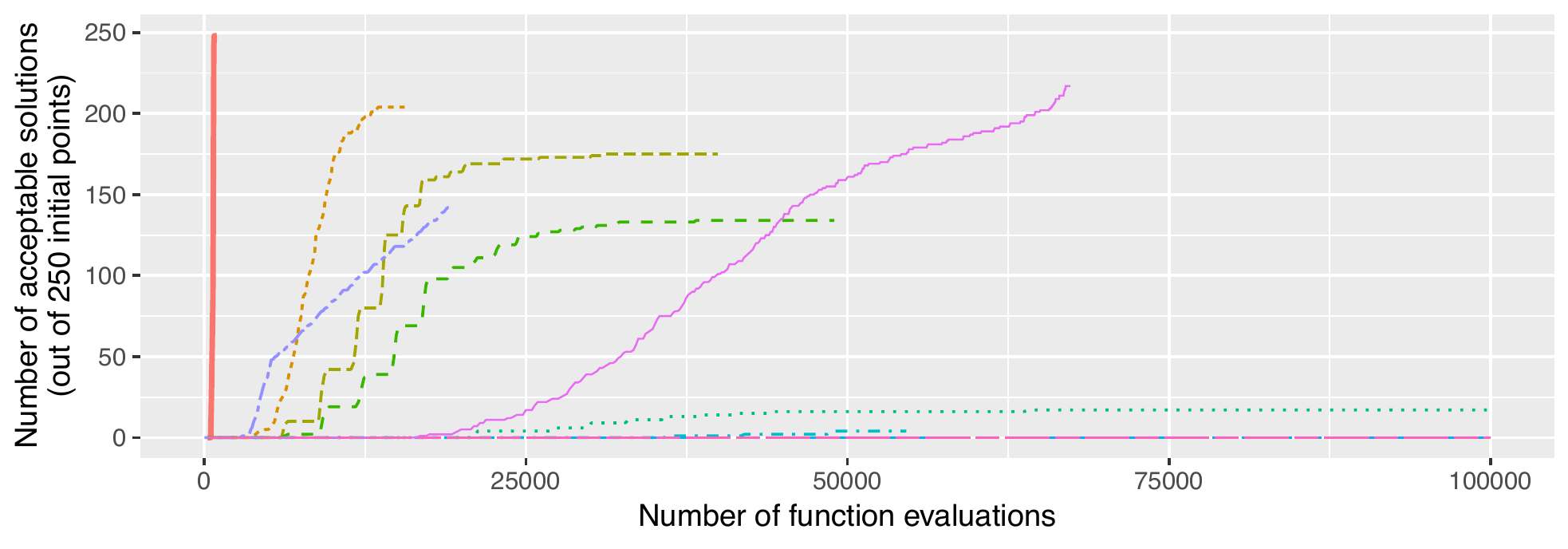}\label{fig::Example2_speed_plot}
\caption{Example 2}
\end{subfigure}
\begin{subfigure}[b]{1\columnwidth}
\includegraphics[width=110mm]{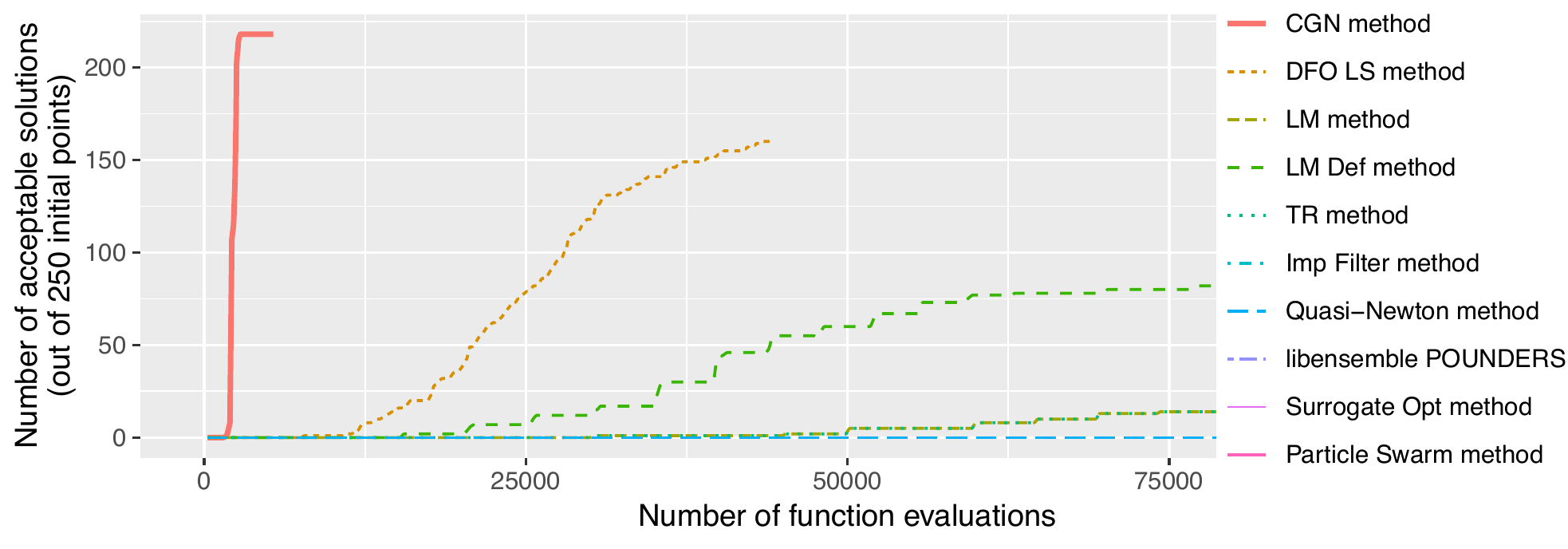}\label{fig::Example3_speed_plot}
\caption{Example 3}
\end{subfigure}
\caption{Number of acceptable minimisers \YA{(out of 250)} found by each method for given total number of function evaluations. Acceptable minimisers are defined by the minimisers with SSR less than the SSR of the parameter that was used to generate the test dataset. (Example1: $\textrm{SSR}<0.0592$, Example2: $\textrm{SSR}<0.0444$, Example3: $\textrm{SSR}<0.0915$)} \label{fig::NumFEAccept}
\end{figure}

 We compared the speed for finding acceptable minimisers.  We define the acceptable minimisers as the minimisers with SSR less than the SSR of the parameter that was used to generate the test dataset.
In Figure~\ref{fig::NumFEAccept}, we show the number of acceptable minimisers found given the total number of function evaluations. As can be seen, CGN is significantly faster than any other method for finding acceptable minimisers.

\begin{figure}

\begin{subfigure}[b]{1\columnwidth}
\includegraphics[width=110mm]{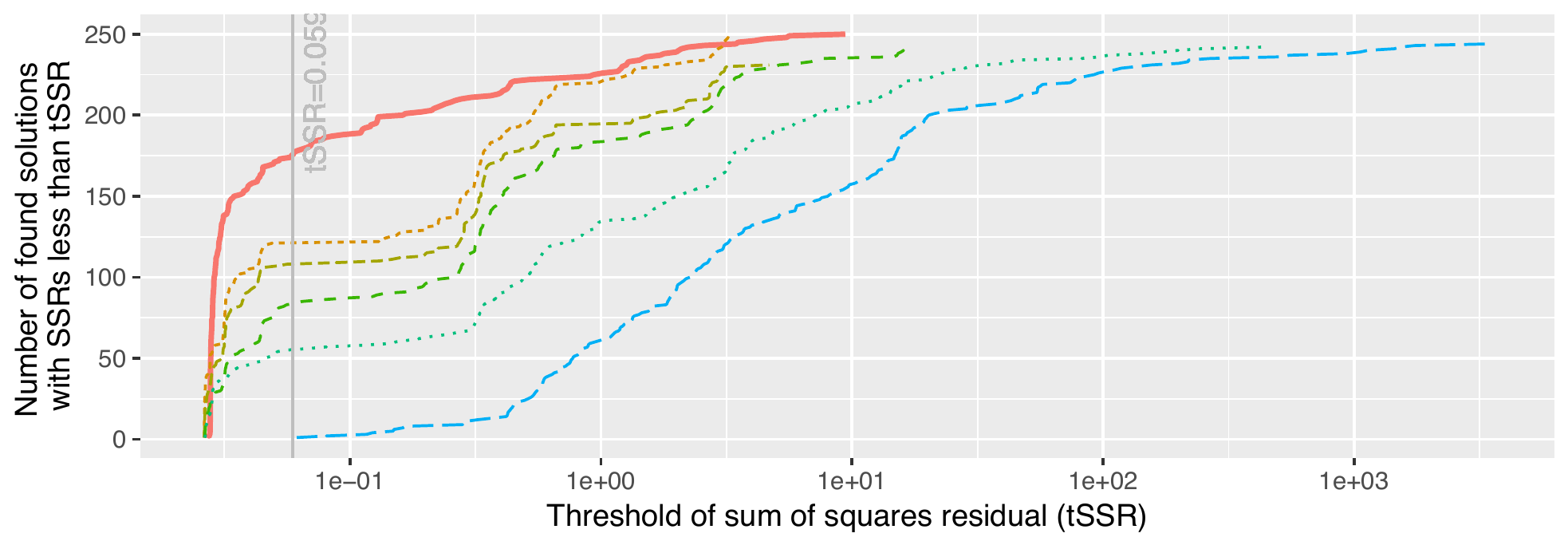}\label{fig::Example1_SSR_plot}
\caption{Example 1}
\end{subfigure}
\begin{subfigure}[b]{1\columnwidth}
\includegraphics[width=110mm]{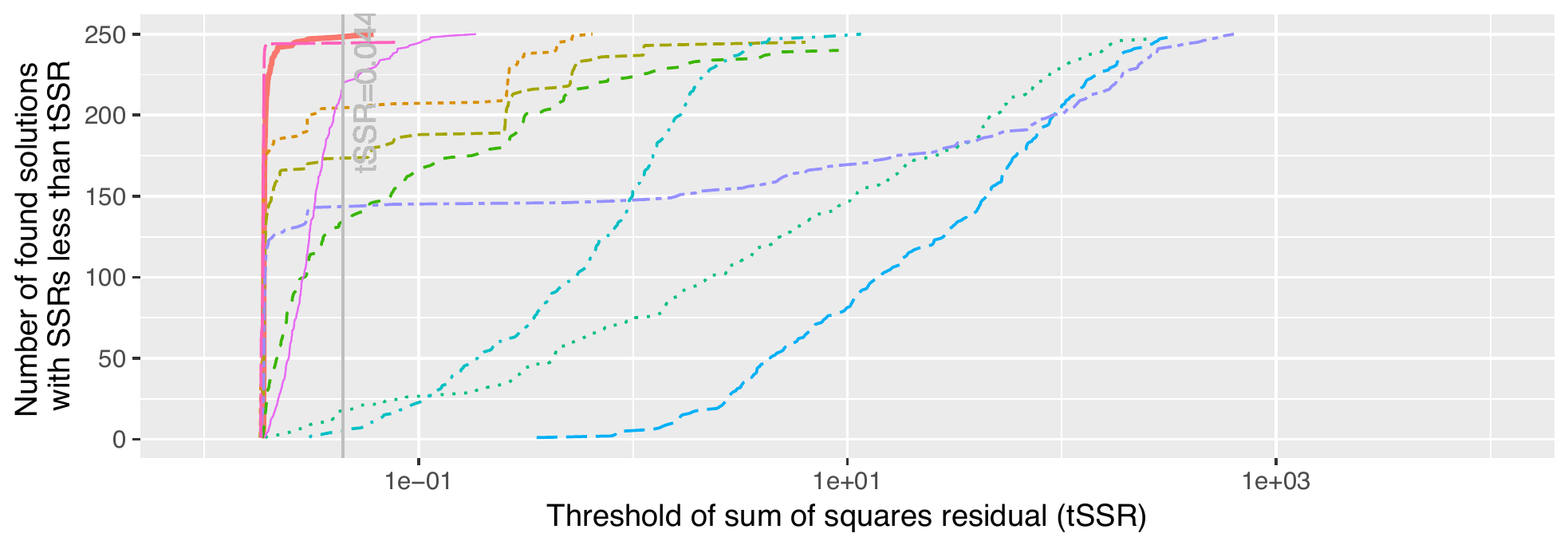}\label{fig::Example2_SSR_plot}
\caption{Example 2}
\end{subfigure}
\begin{subfigure}[b]{1\columnwidth}
\includegraphics[width=110mm]{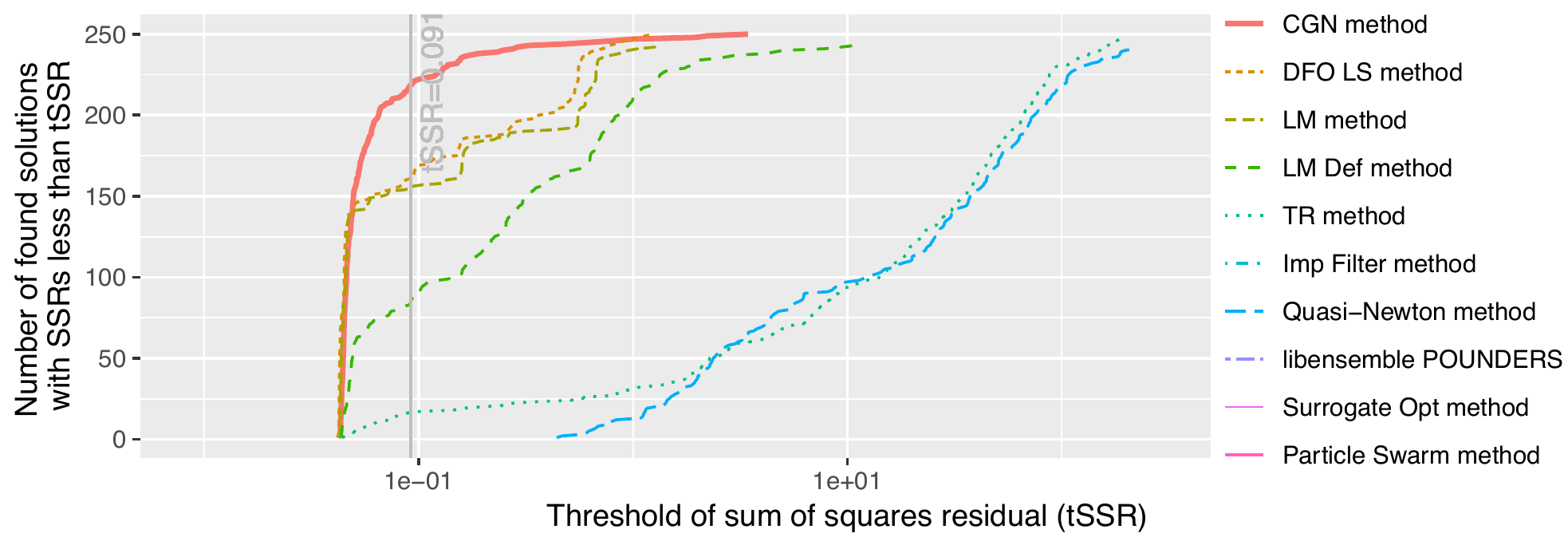}\label{fig::Example3_SSR_plot}
\caption{Example 3}
\end{subfigure}
\caption{Number of solutions \YA{(out of 250)} found by each method for given accuracy threshold (tSSR).  Smaller SSR indicates more accurate solution to the nonlinear-least squares problem.} \label{fig::SSRPlots}
\end{figure}

As this analysis depends on how we define acceptable minimisers, we plotted the number of minimisers found for a given SSR threshold in Figure~\ref{fig::SSRPlots}.  For this analysis each algorithm was run until its default stopping criterion (given in the reference of each method) was met.  We summarised the number of function evaluations that each algorithm required to meet its default stopping criterion in Table~\ref{tab::numFuncEval}.

\YA{As can be seen in Figure~\ref{fig::SSRPlots}, for Examples~1 and 3, CGN finds more approximate minimsers with small SSR than all the conventional algorithms.  For Example~2, as we know that at least one global minimiser is enclosed in the domain defined by $\boldsymbol{x}^\textrm{U}$ and $\boldsymbol{x}^\textrm{L}$, we were able to apply optimisation algorithms with bound constraints (libsensemble with POUNDERS, implicit filter, surrogate optimisation, and particle swarm).  For this problem, particle swarm obtained slightly more approximate minimisers with small SSR than CGN. However, as can be seen in Table~\ref{tab::numFuncEval}, particle swarm required significantly more function evaluations than the CGN (particle swarm: 3,646,900 v.s. CGN: 6,768).  Aside from the particle swarm, CGN outperformed all the other methods in the number of acceptable minimisers found.}

\YA{In addition to above rigorous comparison, we applied two global optimisation algorithms, Genetic algorithm (GA) and DIRECT to Example~2. GA and DIRECT required 120,600 and 100,063 nonlinear function evaluations, respectively, to obtain one solution each.  The SSR of the solutions found by GA and DIRECTwere 0.0197 and 0.244, respectively, while the minimum SSR found by the CGN is 0.0188.  Hence, even to find just one solution, CGN is faster and more accurate than these two methods.}

\begin{table}[ht]
\centering
\begin{tabular}{rrrr}
  \hline
Algorithm & Example 1 &Example 2&Example 3 \\ 
  \hline
 CGN method & 7,782& 6,768& 8,452 \\ 
 DFO LS method & 57,577& 36,804& 70,854 \\ 
 LM method & 72,400& 43,359 & 153,184\\ 
 libensemble POUNDERS & --& 19,068&--\\ 
 LM Def method & 49,822& 43,359& 81,710 \\ 
 Imp Filter method & --& 51,559& -- \\ 
 Trues-Region method & 94,308& 62,136& 103,260 \\ 
 Quasi-Newton method & 124,762& 63,316& 151,774\\ 
 Surrogate Opt method & --& 137,500& -- \\ 
 Particle Swarm method &--& 3,646,900&-- \\
   \hline
\end{tabular}
\caption{Number of total function evaluations.} 
\label{tab::numFuncEval}
\end{table}

%
%


\section{Concluding Remarks}
We proposed the Cluster Gauss-Newton (CGN) \\ 
method, a new derivative free method specifically designed for finding multiple approximate minimisers of a nonlinear least squares problem. 
The development of this algorithm was motivated by the parameter estimation of physiologically based pharmacokinetic (PBPK) models that appears in pharmaceutical drug development.  The particular nature of the model, where the model is over parameterised, and consideration of multiple possible parameters is necessary, motivated us to develop the new method.  The fact that our algorithm obtains multiple \YA{approximate minimiser}s collectively makes it significantly more computationally efficient compared to \YA{existing nonlinear least squares solvers or optimisation algorithms}.  In addition, we observed that in general, CGN obtains minimisers with smaller sum of squared residuals (SSR) than existing algorithms. We demonstrated these advantages using three examples that come from real world drug development projects.

By minimising the assumption on the nonlinear function, where it can be a ``black box'', we have ensured the ease of use and implementation for those who may not have a substantial background in mathematics or scientific computing. We believe this advantage of the proposed method will be appreciated by potential users of the algorithm in industry. In this paper, we used the pharmacokinetics models as examples. However, as we do not assume any particular form of the nonlinear function, we believe the proposed method can be used for many other mathematical models in various scientific fields. 

Software for the CGN method proposed in this paper is available at \url{https://sourceforge.net/p/cluster-gauss-newton-method/code/}.
GUI software is available at \url{http://www.bluetree.me/CGNmethod_for_PBPKmodels}.
Matlab code is available at\\
\url{https://www.mathworks.com/matlabcentral/fileexchange/68798}.

\section{Acknowledgement}
We would like to thank Professor Akihiko Konagaya for giving us the opportunity to initiate this research.  \YA{We would like to thank Dr. Jeffrey Larson for his help on numerical experiment using libensemble. We would like to thank Professor Hiroshi Yabe and Dr. Ryota Kobayashi for useful discussions. We would also like to thank the referees for valuable advice which helped to greatly improve the paper. Ken Hayami was supported in part by JSPS KAKENHI Grant Number 15K04768. Yasunori Aoki is currently employed by AstraZeneca.}

\appendix

\section{Motivating Example}\label{sec::motivatingExample}

\begin{figure}
\centering
\includegraphics[width=120mm]{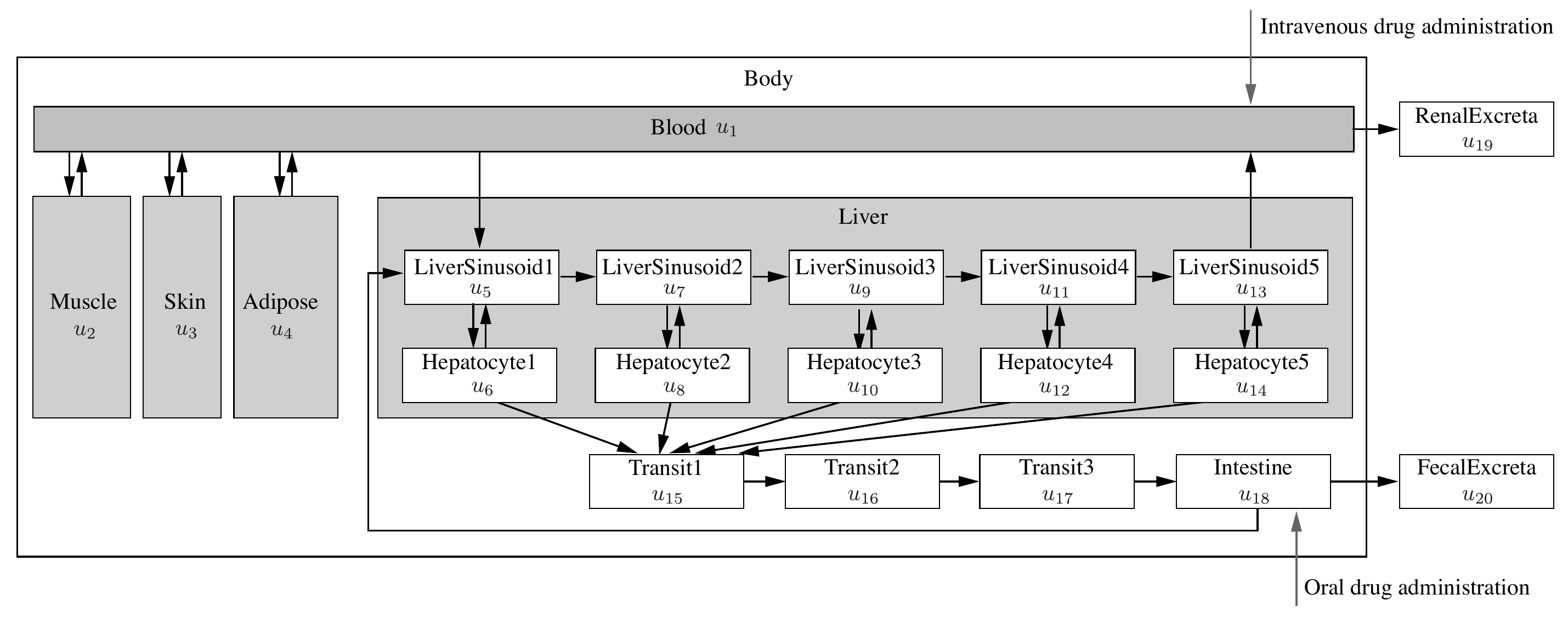}
\caption{A schematic diagram of a physiologically based pharmacokinetic model. (Arrows represent the movement of the drug to a different part of the body. Variables $u_i$ are the drug concentration or the amount of the drug in each compartment. The body is divided into
Blood, Muscle, Skin, Adipose, Liver and Intestine. The liver is further divided into ten compartments to model the complex drug behaviour in the liver. The intestine is divided into four compartments, three transit compartments and one intestine compartment, to model the time it takes for the drug to reach the intestine.)}\label{fig::compartmentDiagram}
\end{figure}

\begin{figure}[h]
\centering
\includegraphics[width=80mm]{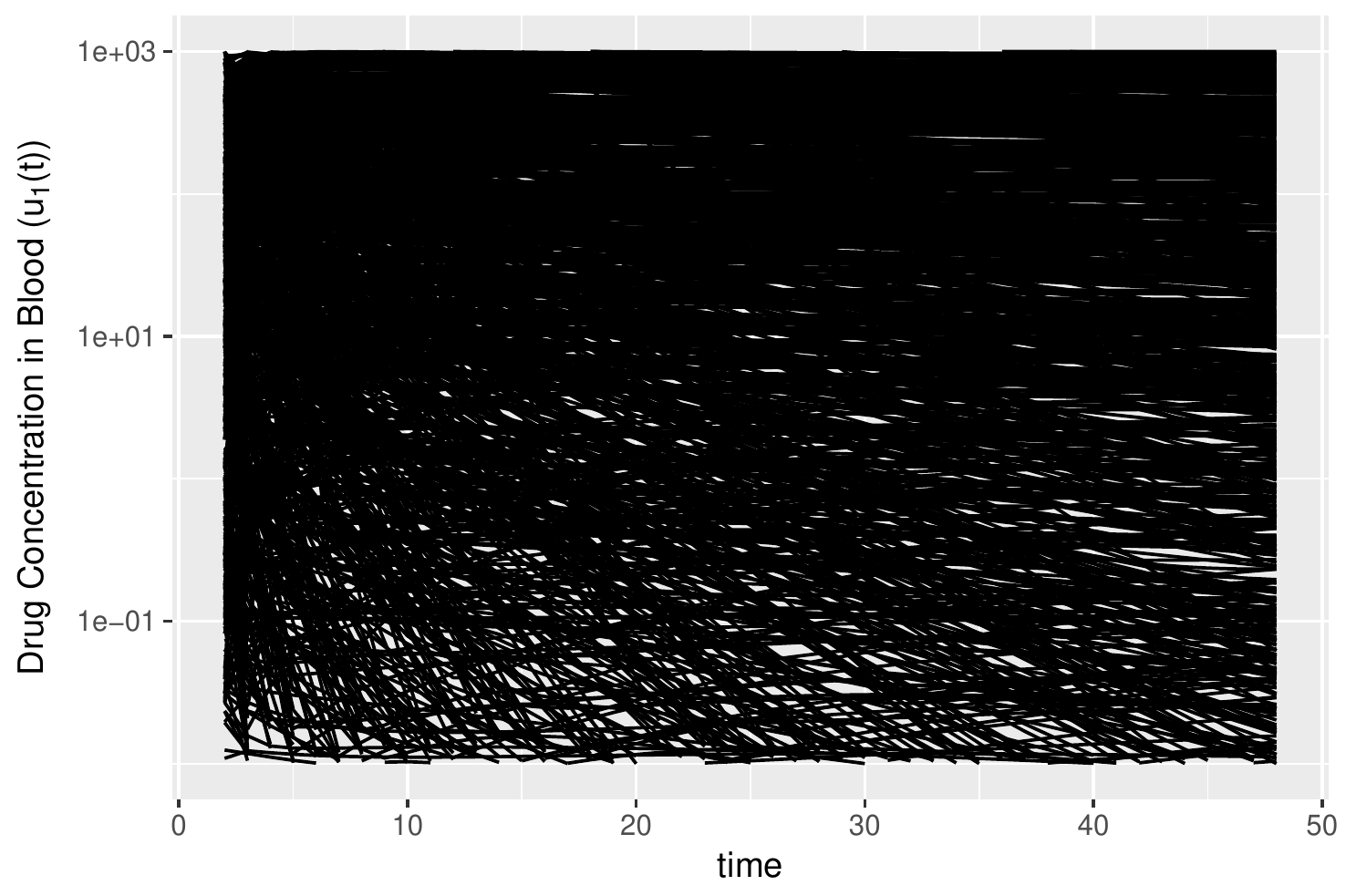}
\caption{Simulation of the drug concentration in the blood plasma using the parameters that are naively sampled from the range of possible kinetic parameters. Note that these simulations do not give any useful information.}\label{fig::simInitPara}
\end{figure}

In this appendix, we introduce a motivating example to illustrate how the proposed method could be used in practice. We consider a scenario where a newly developed drug is tested for the first time in a human. Before the drug is given to a human, the biochemical properties of the drug are studied in-test-tube (in-vitro) and in-animal experiments. However, how the drug behaves in the human body is still uncertain. Based on the results of in-vitro and in-animal experiments, the team decides that 100mg is a safe amount of the drug to be given to a human and the experiment is conducted with a healthy normal volunteer, and the drug concentration in the blood plasma is measured at points in time. Using these measurements, we estimate the multiple possible model parameters which can be used to simulate various scenarios. The following workflow can be envisioned:

\textbf{1: Construct a mathematical model based on the understanding of the physiology and biochemical properties of the drug.}\\
 In this example, we use the model presented in \cite{Watanabe2009}. The mathematical model is depicted in Figure~\ref{fig::compartmentDiagram}, and it can be written as a system of nonlinear ordinary differential equations with 20 variables.
There are two types of model parameters in this model: physiological parameters and kinetic parameters. Examples of the physiological parameters are the sizes of the organs or the blood flow rates between the
organs. As the human physiology is well studied and these parameters usually do not depend on the drugs, we can assume these parameters to be known. The kinetic parameters, such as, how fast the drug gets excreted from the body or how easily it binds to tissues are the parameters that depend on the drug and usually are not very well known. Before the first in-human experiment, the drug development team characterises these parameters using an organ in test-tubes or by administering the drug to an animal. However, these parameters can differ from animal to human, so we do not have a very accurate estimate of these parameters. The differences in these parameters between a human and an animal can be several orders of magnitude. Figure~\ref{fig::simInitPara} depicts plots of the drug concentration simulation where the kinetic parameters are sampled within a reasonable range of the parameters. As can be seen in Figure~\ref{fig::simInitPara}, we cannot obtain any useful information just by randomly sampling the kinetic parameters from the feasible range.

\begin{figure}
\centering
\includegraphics[width=120mm]{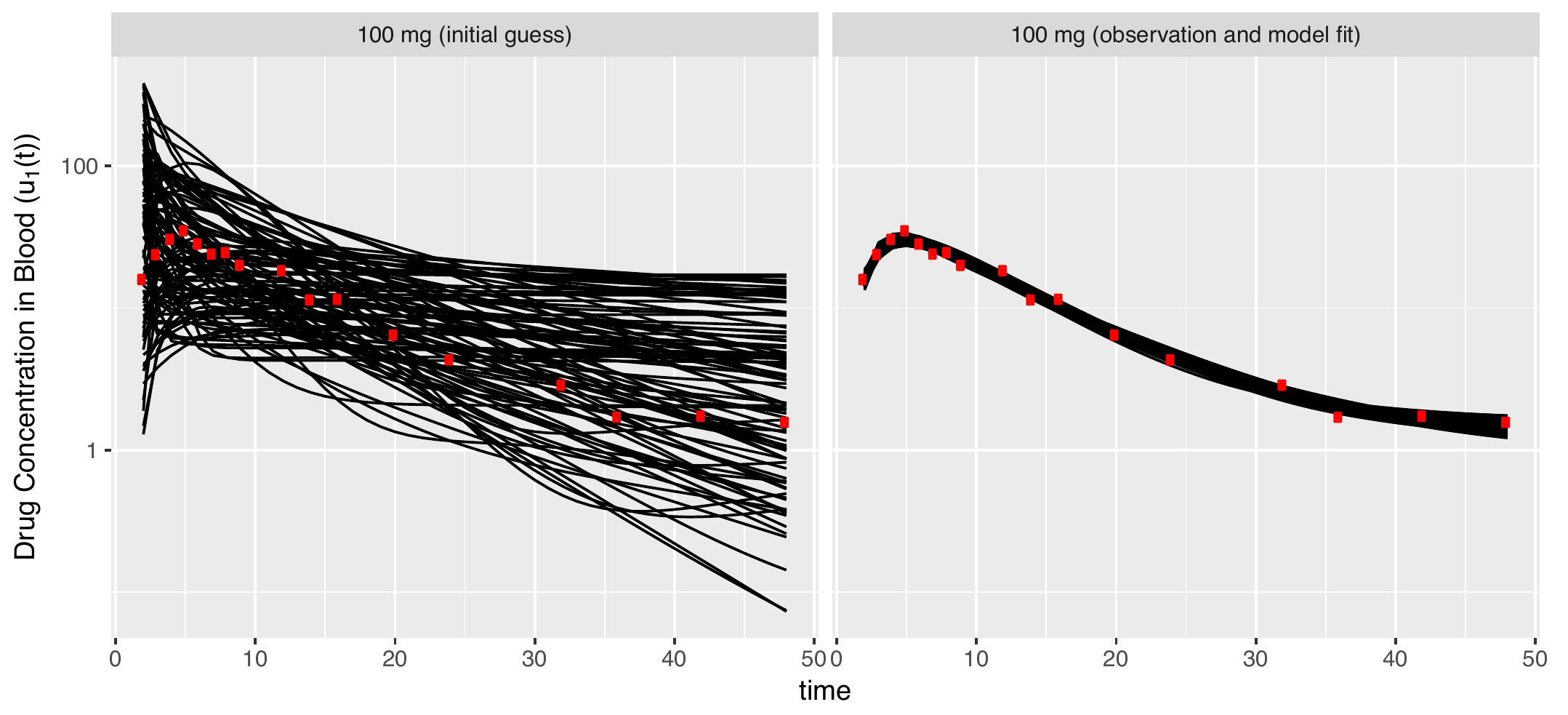}
\caption{
Plot of the simulation of drug concentration (black solid line) with observations (red dot). 
Simulation results are based on the parameters for the initial iterate for CGN and the parameters found after 20 iterations of CGN.
 In the left panel, the simulation results based on the top 100 sets of the parameters (out of 1000 parameter sets in the cluster) from the initial cluster are shown. In the right panel, the simulations results based on the top 100 sets of the parameters (out of 1000 parameter sets in the cluster) after 20 iterations of CGN are shown. Hence, top 100 means the 100 smallest sum of squared residuals (SSR) between the simulation and observed data.}\label{fig::simulation}
\end{figure}

\begin{figure}
\centering
\includegraphics[width=120mm]{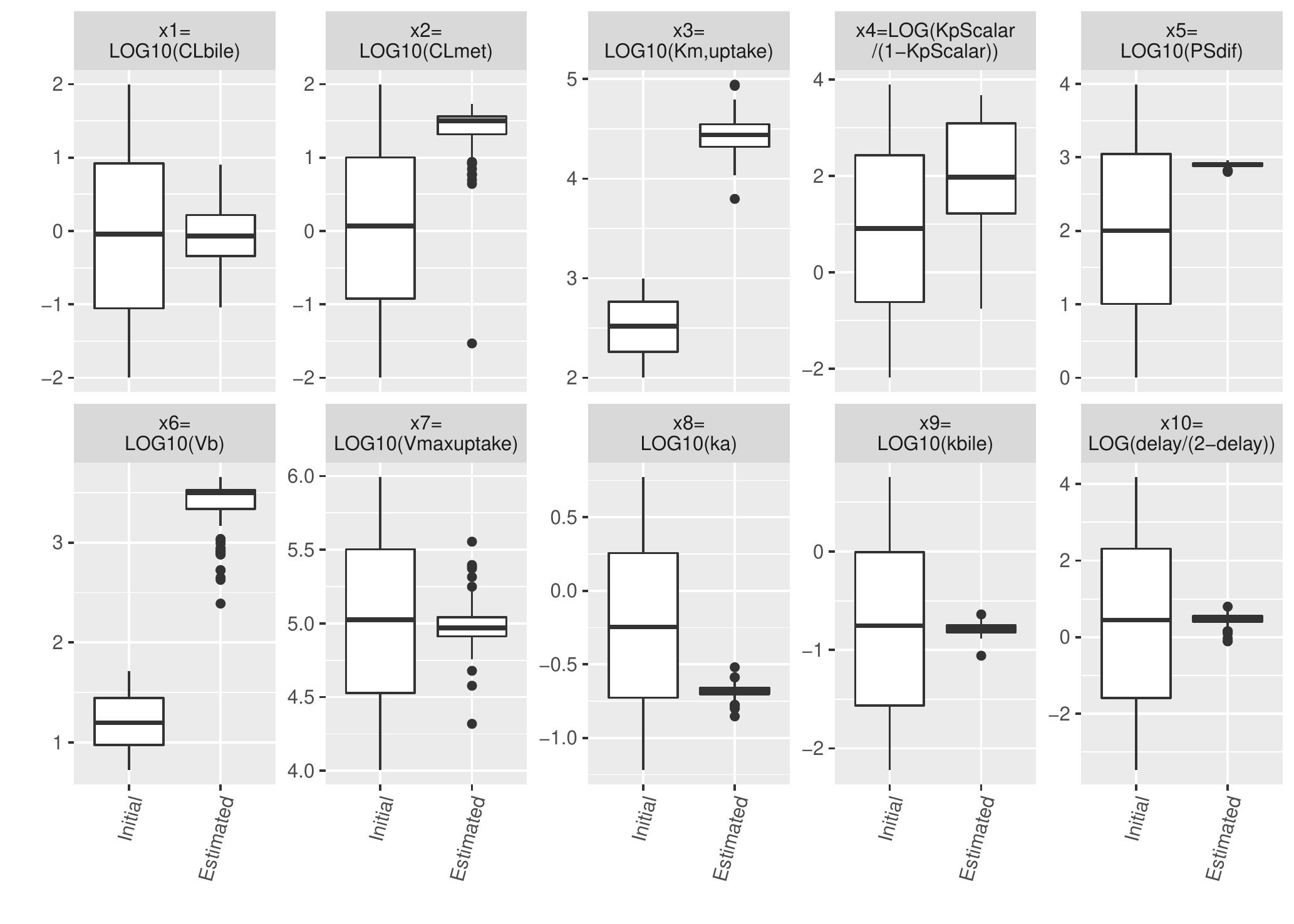}
\caption{Box plots of top 100 parameters (the same parameter sets as the one used to plot Figure~\ref{fig::simulation}) from the initial cluster (left) and the cluster after 20 iterations of CGN (right). Note the distributions of $x5, x8, x9, x10$ clearly shrunk after 20 iterations while the distribution of $x4$ did not change noticeably. (Box plot: The edges of the boxes are the 75th and 25th percentiles. The line in the box is the median, and the whiskers extend to the largest and the smallest value within the 1.5 times the inter-quartile range. Dots are the outliers outside the whiskers.)}\label{fig::boxplot}
\end{figure}

\begin{figure}
\centering
\includegraphics[width=120mm]{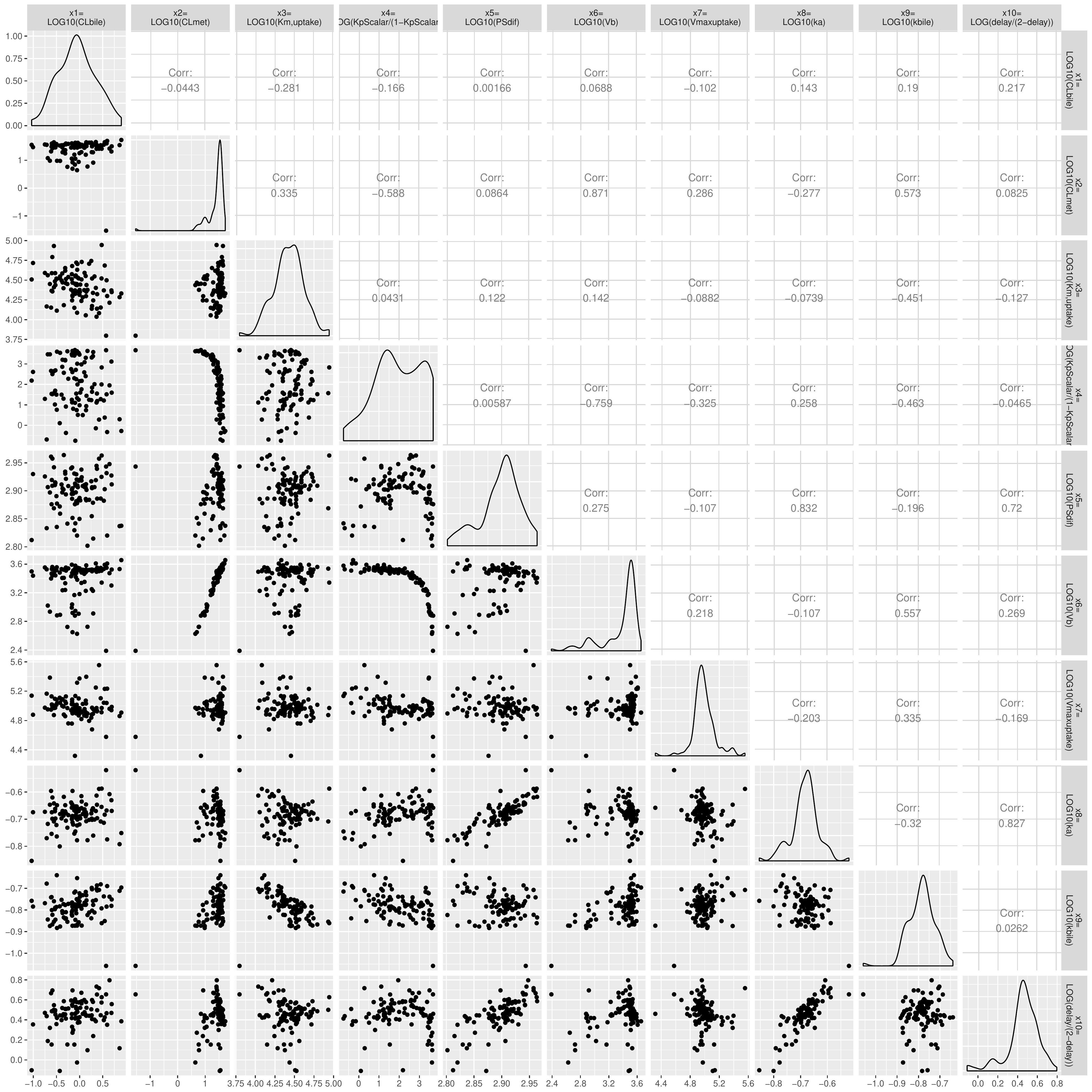}
\caption{Scatter plots of parameters. As can be seen in this example, we can find parameter-parameter correlations of some of the parameters found by CGN. (parameters whose corresponding SSR is the least 100 within the cluster of 1,000, i.e., the same parameter sets as the one used to plot the right hand side of Figure~\ref{fig::simulation}.)}\label{fig::scatterPlots}
\end{figure}

\textbf{2: Sample multiple possible model parameter sets that fit the model prediction of the drug concentration to the observed data from the 100mg experiment.} We now use the observed data from the experiment where 100mg of the drug was given to a human. The red dots in Figure~\ref{fig::simulation} depict the observed data. The left panel of Figure~\ref{fig::simulation} shows some of the simulation results using the parameter sets of the initial iterate of CGN. The right panel of Figure~\ref{fig::simulation} shows some of the simulation results after 20 iterations of CGN. As can be seen in Figure~\ref{fig::simulation}, CGN can find multiple sets of parameters that fit the observed data. The parameter values are depicted in the box plots in Figure~\ref{fig::boxplot}. As can be seen in Figure~\ref{fig::boxplot}, after 20 iterations of CGN, the distribution of some of the parameters shrinks significantly suggesting these parameters can be identified from the observations while the distribution of some of the parameters are unchanged indicating that these parameters cannot be identified from the observation. In Figure~\ref{fig::scatterPlots}, we show scatter plots of the parameters found by CGN. As can be seen in Figure~\ref{fig::scatterPlots}, even if the parameter cannot be identified from the observation, some nonlinear relationships can occasionally be identified between the parameters.

\section{Relation with previous work}\label{header-c1650}

We initially proposed the Cluster Newton (CN) method for sampling multiple solutions of a system of nonlinear equations where the number of unknowns is more than the number of observations \cite{Aoki2011, Aoki2014, Gaudreau2015}.  

The parameter fitting problems presented in \cite{Yoshida2013, Toshimoto2017, Kim2017} uses CN  to fit complex PBPK models to the drug (and its metabolite) concentration measurements over time. As the drug concentration was measured repeatedly from the patient, there is a larger number of observations than the number of unknown parameters. CN assumes an underdetermined problem, where the number of observations is less than the number of parameters. To use CN, in \cite{Yoshida2013, Toshimoto2017, Kim2017}, the authors constructed a summary value called Area Under the time-concentration Curve (AUC) to reduce the number of observations and fit the model to the summary value using CN. The AUC is essentially the time integral of the drug concentration from the time drug is administered to infinity. After finding multiple possible parameters that match the AUC of the model and the observations, Yoshida et al. and Toshimoto et al. \cite{Yoshida2013, Toshimoto2017} selected the parameter sets that fit the time-course drug concentrations reasonably well from these parameter sets found using CN. Kim et al. \cite{Kim2017} used a parameter set obtained using CN as the initial iterate for the Levenberg-Marquardt method and fitted the model to the time-course drug concentration data.

Based on the current use of CN, we identified two bottlenecks of the workflow employed in \cite{Yoshida2013, Toshimoto2017, Kim2017}. Firstly, CN solves a system of underdetermined system of nonlinear {\bf equations}. Hence, the model needs to be constructed in a way that the observation and the model prediction match {\bf exactly}. Hence, if there is a model misspecification or significant measurement error that influence the summary values (e.g. AUC) sufficiently so that there is no model prediction that exactly matches the observation derived summary value, then CN breaks down.

Secondly, in order to obtain the parameter sets that reasonably fit the original data (e.g., time-course concentration data), we need to obtain many parameter sets that fit the summary values (e.g., AUC). This is simply because we need to randomly sample from (number of parameters)-(number of summary values) dimension space to obtain the desired parameter sets. As a result, \cite{Toshimoto2017} had to find 500 000 parameter sets that fit the summary value (AUC) using CN and then was able to obtain 30 parameter sets that reasonably fit the original data (time-course drug concentration).

\begin{figure}
\includegraphics[width=70mm]{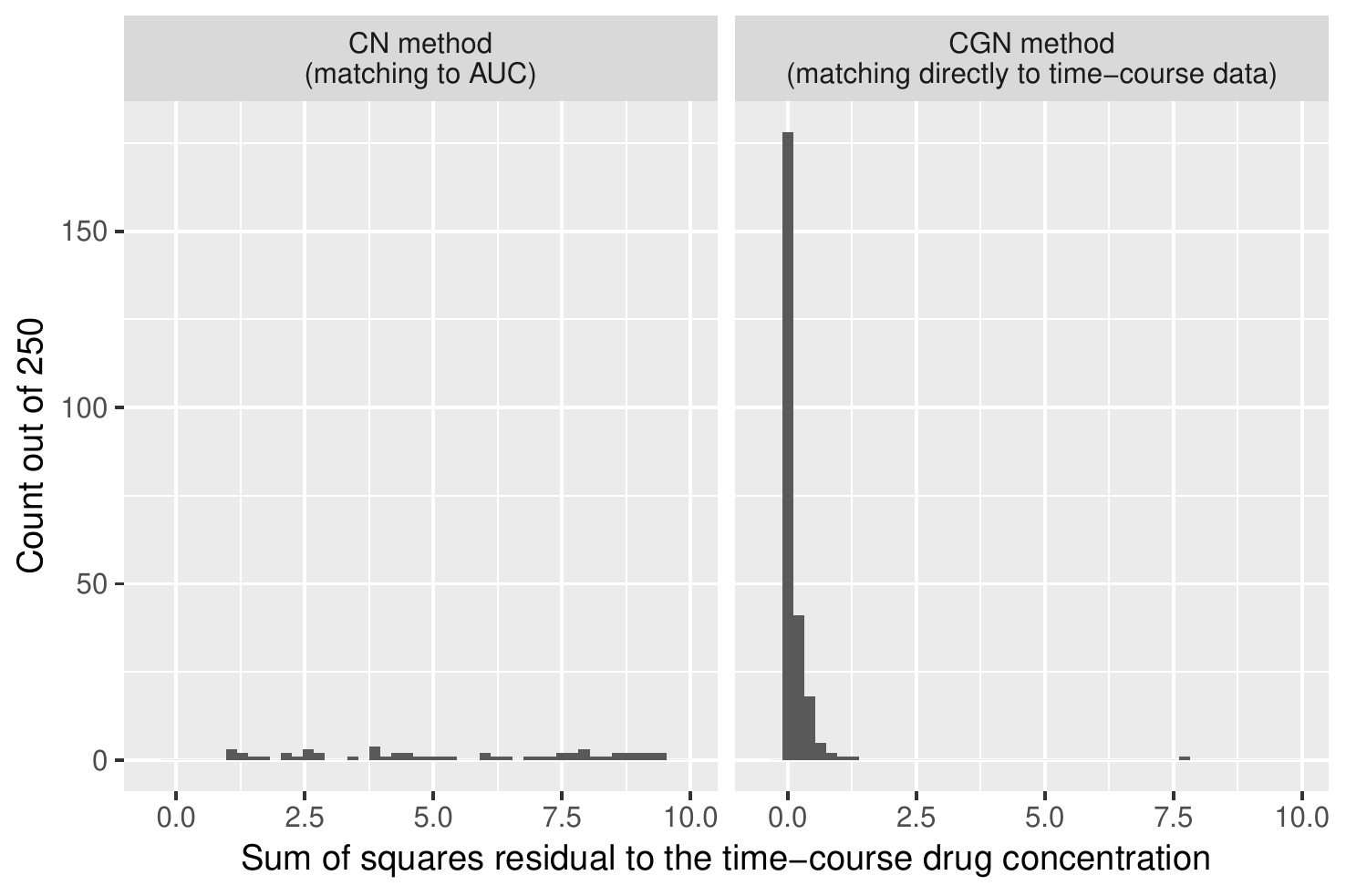}	
\caption{Comparison of the accuracy of the model fit of the Motivating Example (cf. Appendix~\ref{sec::motivatingExample}) using the approach employed in \cite{Yoshida2013, Toshimoto2017} using the Cluster Newton (CN) method and the proposed Cluster Gauss-Newton (CGN) method.}\label{fig:comparison_new_old_cnmethod}
\end{figure}

To overcome these bottlenecks, we suggested these authors to formulate these parameter estimation problems as nonlinear least squares problem and we created the Cluster Gauss-Newton method (CGN). The new method efficiently obtains multiple possible parameters by solving a nonlinear least squares problem so that the method does not break down even if the model does not exactly match the observation. Also, the new algorithm does not require the number of observations to be less than the number of parameters so that there is no need to summarise the observation and can directly use the original observations (e.g., time-course concentration measurements).  As can be seen in Figure~\ref{fig:comparison_new_old_cnmethod}, the new CGN algorithm enables us to sample multiple possible parameters significantly more accurately than CN.

\section{Numerical experiments on the influence of the weight of the linear approximation} \label{app::NumExpWeight}

\begin{figure}
\centering
\includegraphics[width=120mm]{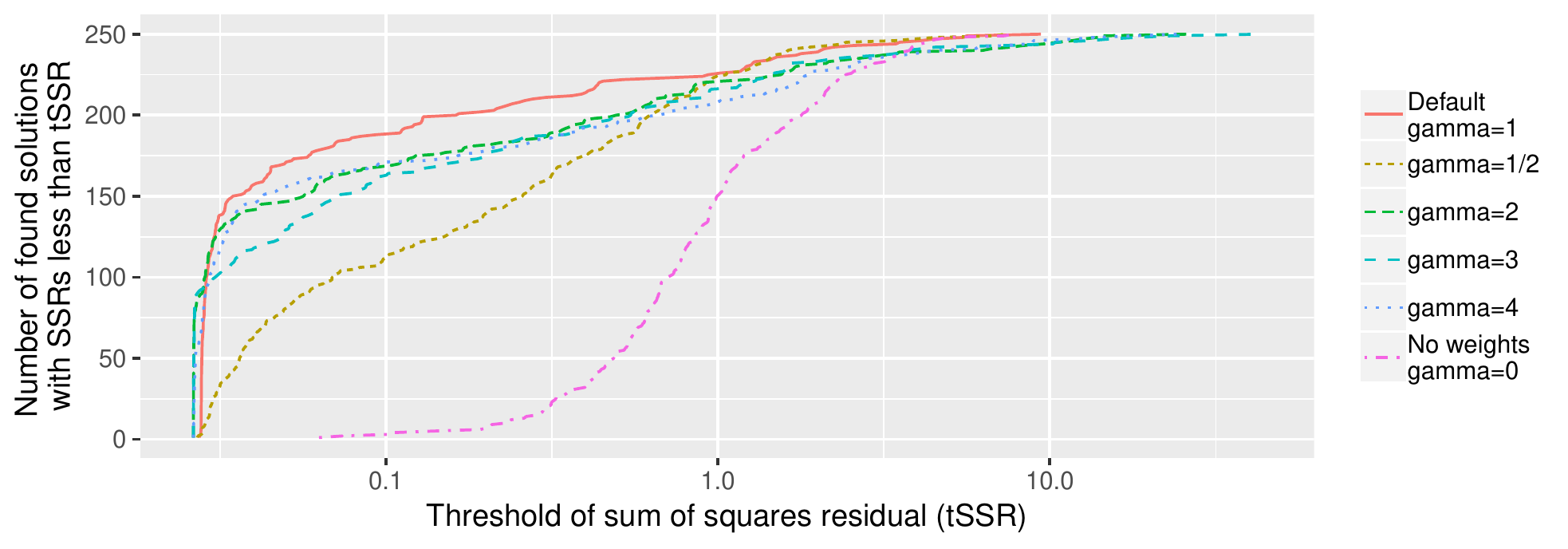}

\caption{Number of solutions found for given accuracy threshold (tSSR) using various weights for the linear approximation (i.e., various $\gamma$).  Smaller SSR indicates more accurate solution to the nonlinear-least squares problem.} \label{fig::sortedSSR_varioutGamma}
\end{figure}

To illustrate the influence of the weight $d_{j(i)}^{(k)}$ of the weighted linear approximation in 2-1) of the algorithm (equations (\ref{eq::linAppSum}) and (\ref{eq::weights}) ), we conducted numerical experiments using Example 1 of section \ref{header-c1748}.  For this numerical experiment, we varied the parameter $\gamma \geq 0$ in Equation~\eqref{eq::weights}:
\begin{eqnarray}d_{j(i)}^{(k)}=\left\{ 
\begin{array}{ll}\left(\frac{1}{
\sum_{l=1}^{n}(({x}_{lj}^{(k)}-{x}_{li}^{(k)})/(x_l^\textrm{U}-x_l^\textrm{L}))^2}\right)^\gamma & \textrm{if }j\neq i\\
0&\textrm{if }j = i
\end{array}
\right. .\end{eqnarray}
 In Figure~\ref{fig::sortedSSR_varioutGamma}, we show the number of solutions found by CGN for given accuracy threshold.  
 As can be seen from the figure
 , the weight for the linear approximation improves the accuracy and the speed of CGN.  Note that $\gamma =0$, which corresponds to giving equal weights to all the cluster points is not optimal.  In this example, $\gamma\geq 1$ gave good results.

\section{ODE expressions of the PBPK model used as the motivating example}\label{app::PBPKmodel}
In this subsection we explicitly write out the mathematical model used in the motivating example in Appendix~\ref{sec::motivatingExample} as a system of ODEs.  The other PBPK models used in this paper can be similarly written as a system of ODEs.  For the model parameters appearing in these ODEs, we follow the usual notations used in pharmacokinetics and we will provide brief descriptions of these parameters.  For more detailed description of these parameters, we refer the reader to the standard textbooks in pharmacokinetics (e.g., \cite{rowland2011clinical}).\\ \\
\noindent {\bf ODE for the blood compartment}
\begin{eqnarray}
\dee{u_1}{t} &=& \frac{Q_h (u_{13} -  u_1) - CL_r u_1 - Q_m (u_1 - \frac{u_2}{
Kp_m Kp_{scalar}}) - Q_s (u_1 - \frac{u_3}{Kp_s Kp_{scalar}}) }{V_b}\nonumber\\
&&- \frac{Q_a (u_1 - \frac{u_4}{Kp_a Kp_{scalar}})}{V_b}\nonumber
\end{eqnarray}

$Q_h$ is the blood flow rate through the liver. The unit is typically L/hr.

$CL_r$ is the rate of the drug being excreted as urine. The unit is typically L/hr.

$V_b$ is the volume of the blood compartment. The unit is typically L.

$Q_m, Q_s, Q_a$ are the blood flow rates to muscle, skin, and adipose, respectively. The unit is typically L/hr.

$V_m, V_s, V_a$ are the volumes of muscle, skin, and adipose compartments, respectively. The unit is typically L.

$Kp_m, Kp_s, Kp_a$ are the partition coefficients of the drug for the muscle, skin and adipose, respectively.  The partition coefficient is the ratio between the drug concentration in the blood plasma and the tissue at equilibrium when the tissue is in contact with the blood.  We usually assume that these values can be measured by \textit{in-vitro} experiments. Hence, they are usually fixed parameters in the model.  

$Kp_\textrm{scalar}$ is a scaling factor for the partitioning. \\ \\
{\bf ODE for the muscle compartment}
\begin{eqnarray}
\dee{u_2}{t}= \frac{Q_m}{V_m}  \left(u_1 -  \frac{1}{ Kp_m  Kp_\textrm{scalar}}u_2\right)
\end{eqnarray}
{\bf ODE for the skin compartment}
\begin{eqnarray}
\dee{u_3}{t}=\frac{Q_s}{V_s} \left(u_1 - \frac{1}{Kp_s Kp_\textrm{scalar}}u_3 \right)
\end{eqnarray}
{\bf ODE for the adipose (fat) compartment}
\begin{eqnarray}
\dee{u_4}{t}=\frac{Q_a}{V_a} \left(u_1 - \frac{1}{Kp_a Kp_\textrm{scalar}}u_4 \right)
\end{eqnarray}
{\bf ODEs for the liver compartments}\\
For a compartment representing the blood vessel in the liver (Liver Sinusoid 1)
\begin{eqnarray}
\dee{u_5}{t}=-\frac{\frac{\textrm{Vmax}_\textrm{uptake}}{Km_\textrm{uptake} + u_5} + f_b \, PS_\textrm{dif} }{V_{hc}}u_5 
+\frac{ f_h \, PS_\textrm{dif} }{V_{hc}} u_6 
+\frac{Q_h (u_1 - u_5)  + ka \, u_{18}}{\frac{V_{hc}}{5}}
\end{eqnarray}
$ka$ is the rate of the drug being absorbed from the intestines and transported to the liver through the portal vein. The unit is typically L/hr.

$V_{hc}$ is the volume of the blood in the blood vessel.

$PS_\textrm{diff}$ is the diffusion constant between the liver sinusoid compartment and the hepatocyte compartment.

$f_b$ is the ratio of the drug that is not bound to the protein in blood in the blood vessel in the liver.

$f_h$ is the ratio of the drug that is not bound to the protein in the liver cells.  Only the portion of the drug that is not bound to the protein can permeate in and out of the liver cells.\\

For compartments representing the blood vessel in the liver (Liver Sinusoid 2, 3, 4, 5)
\begin{eqnarray}
\dee{u_{7+2i}}{t}=- \frac{ \frac{\textrm{Vmax}_\textrm{uptake} }{Km_\textrm{uptake} + u_{7+2i}} + f_b \, PS_\textrm{dif} }{V_{hc}}u_{7+2i}
+ \frac{ f_h \, PS_\textrm{dif}}{V_{hc}}u_{8+2i}+\frac{Q_h  (u_{5+2i} -  u_{7+2i})}{\frac{V_{hc}}{5}}\nonumber
\end{eqnarray}
for $i=0,1,2,3$.

$\textrm{Vmax}_\textrm{uptake}$ and $Km_\textrm{uptake}$ are constants for the Michaelis-Menten kinetics where the active uptake of the drug from the blood to the liver cells by membrane transport protein (also known as ``transporter") is modelled.\\

Compartments representing liver cells (Hepatocyte 1, 2, 3, 4, 5)
\begin{eqnarray}
\nonumber
\dee{u_{6+2i}}{t}=\frac{\frac{\textrm{Vmax}_\textrm{uptake} }{ Km_\textrm{uptake} + u_{5+2i}} + f_b \, PS_\textrm{dif} }{V_{he} }u_{5+2i}
-\frac{f_h ( PS_\textrm{dif} + (\textrm{CL}_\textrm{met} +\textrm{CL}_\textrm{bile}) )}{V_{he} }u_{6+2i}
\end{eqnarray}
for $i=0,1,2,3,4$.

$V_{he}$ is the volume of liver cells.\\ \\
{\bf ODE for the transit compartment 1}
\begin{eqnarray}
\dee{u_{15}}{t}= f_h \, \textrm{CL}_\textrm{bile} \, \frac{u_6 + u_8 + u_{10} + u_{12} + u_{14}}{5} - k_\textrm{bile}\, u_{15}
\end{eqnarray}
$u_{15}$ is the amount of the drug in the transit compartment 1.  The transit compartment can be thought of as a fictitious compartment that introduces the delay of the drug delivery to to the intestines.

$\textrm{CL}_\textrm{met}$ is a clearance (speed that drug is cleared out) due to the metabolisation of the drug of the subject. The unit is typically L/hr. 

$\textrm{CL}_\textrm{bile}$ is a clearance due to the biliary excretion of the drug of the subject. The unit is typically L/hr. 

$k_\textrm{bile}$ is a diffusion coefficient between the transit compartments. The unit is typically 1/hr.  

$f_h$ is the fraction of the drug cleared by liver, it is typically unit-less.\\
\\
{\bf ODE for the transit compartment 2}
\begin{eqnarray}
\dee{u_{16}}{t}=k_\textrm{bile} \, (u_{15} - u_{16})
\end{eqnarray}

$u_{16}$ is the amount of the drug in the transit compartment 2.

$k_\textrm{bile}$ is a diffusion constant where the drug is diffused into bile, the unit is typically 1/hr.\\ \\
{\bf ODE for the transit compartment 3}
\begin{eqnarray}
\dee{u_{17}}{t}=k_\textrm{bile} \, (u_{16} - u_{17})
\end{eqnarray}

$u_{17}$ is the amount of the drug in the transit compartment 3.\\ \\
{\bf ODE for the intestine compartment}
\begin{eqnarray}
\dee{u_{18}}{t}=k_\textrm{bile} \, u_{17} - \frac{ka }{FaFg} \, u_{18}
\end{eqnarray}

$u_{18}$ is the amount of the drug in the intestine.

$FaFg$ is the fraction of the drug that gets absorbed by the intestines.

See also Appendix \ref{app::Repara}

\section{Re-parameterisation}\label{app::Repara}

For the mathematical model we introduced in Appendix~\ref{app::PBPKmodel}, some of the parameters may be known for the drug of interest.  For the motivating example presented in Appendix~\ref{sec::motivatingExample}, we assumed the following parameters to be known:

\begin{eqnarray}
CL_r=0,\qquad
FaFg=0.55,\qquad
Kp_a=0.086,\qquad
Kp_m=0.113\\
Kp_s=0.478,\qquad
Q_a=15.61,\qquad
Q_h=86.94,\qquad
Q_m=44.94\\
Q_s=17.99,\qquad
V_a=10.01,\qquad
V_{hc}=1.218,\qquad
V_{he}=0.469\\
V_m=30.03,\qquad
V_s=7.77,\qquad
f_b=0.00617,\qquad
f_h=0.012
\end{eqnarray}

Also, most of the parameters have known physiologically or chemically possible range (e.g., clearances cannot be negative).  Thus, for those parameters, we re-parameterise to impose known constraints. For the motivating example presented in Appendix~\ref{sec::motivatingExample}, we re-parameterised as in the following:
\begin{eqnarray}
	\textrm{CL}_\textrm{bile}=10^{x_1},\qquad
\textrm{CL}_\textrm{met}=10^{x_2},\qquad
Km_{\textrm{uptake}}=10^{x_3},\\
Kp_\textrm{Scalar}=\frac{e^{x_4}}{1+e^{x4}}\qquad
PSdif=10^{x_5},\qquad
Vb=10^{x_6},\\
Vmax_{\textrm{uptake}}=10^{x_7},\qquad
ka=10^{x_8}\qquad
k_\textrm{bile}=10^{x_9}.
\end{eqnarray}

\section{Numerical experiments on the influence of the regularisation} \label{app::NumExpReg}
\begin{figure}
\centering
\includegraphics[width=120mm]{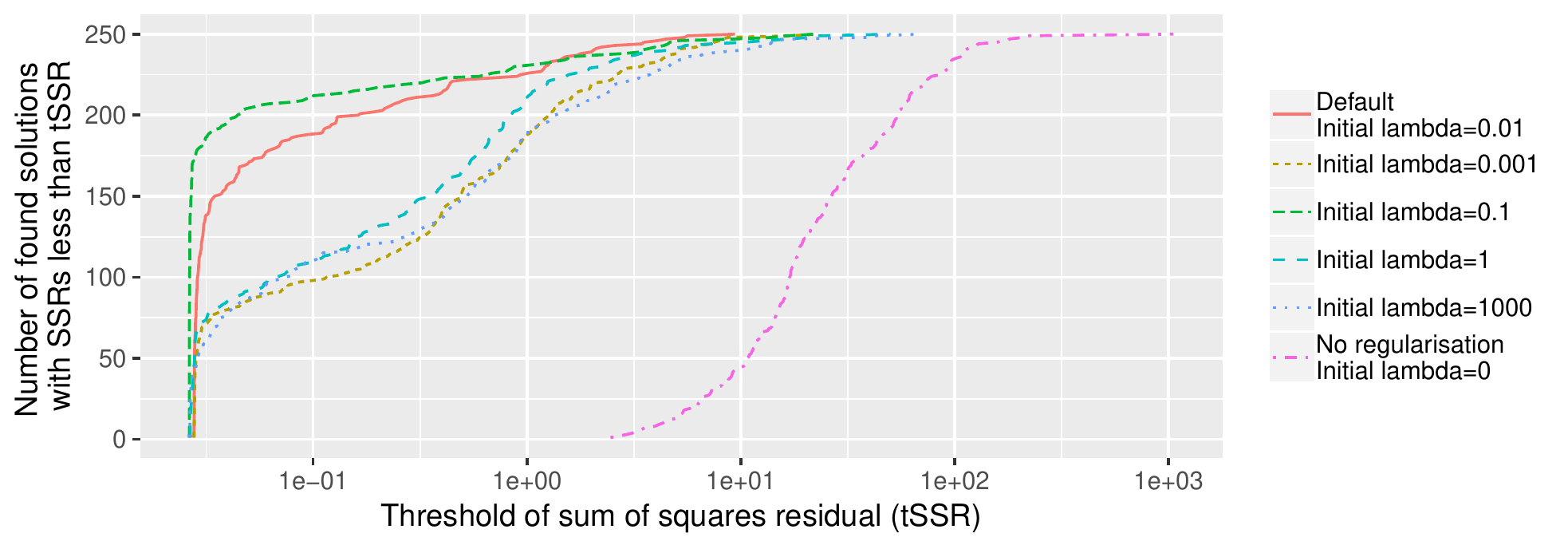}

\caption{Number of solutions found for given accuracy threshold (tSSR) using various initial lambda ($\lambda_\textrm{init}$).  Smaller SSR indicates more accurate solution to the nonlinear-least squares problem.} \label{fig::sortedSSR_varioutLambda}
\end{figure}

To illustrate the necessity and the influence of the regularisation in 2-2) of the algorithm ( equation (\ref{eq::regularise}) ), we  conducted numerical experiments using Example 1 of section \ref{header-c1748}. We varied the initial value of the regularisation coefficient $\lambda_\textrm{init}$ and tested CGN.
Figure~\ref{fig::sortedSSR_varioutLambda} shows the number of solutions found by CGN for given accuracy threshold.  
As can be seen from the figure, 
regularisation is necessary for CGN to perform well.  For this example, $\lambda_\textrm{init}=0.1$ gave a good result.

\section{Application to discontinuous nonlinear functions}
The Jacobian based local optimisation methods cannot solve nonlinear least squares problems if the nonlinear function is discontinuous.  On the other hand, the Cluster Gauss-Newton (CGN) method uses the linear approximation of the nonlinear function to capture the global behaviour of the function. Hence, it does not require the nonlinear function to be continuous.

For the following numerical experiment, we consider the case where the nonlinear function is not continuous.  This example is constructed to show that CGN can solve nonlinear least squares problems that the conventional Jacobian based method cannot solve.  We create such a nonlinear function by rounding the nonlinear function of Example 1  (\eqref{eq::Example1_ODE}-\eqref{eq::Example1_IC}) to the first decimal place.

\begin{figure}[h]
\includegraphics[width=120mm]{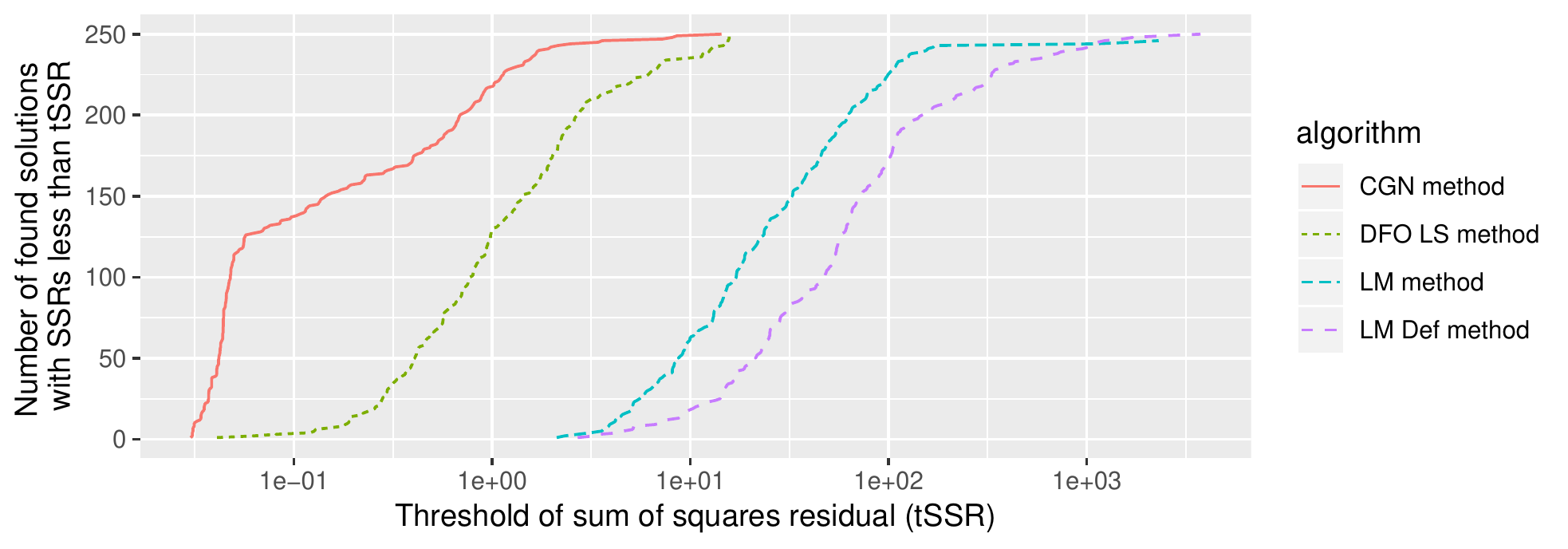}
\caption{Example 1 when the nonlinear function is rounded to one decimal place: 
Number of solutions \YA{(out of 250)} found by the various methods for given accuracy threshold (tSSR).  Smaller SSR indicates more accurate solution to the nonlinear-least squares problem.}\label{fig::Example1SSR_round}
\end{figure}

As can be seen in Figure~\ref{fig::Example1SSR_round}, the CGN method was able to find many accurate solutions.  On the other hand, the LM failed to find any reasonable solution. Unlike derivative based LM, a derivative-free method (DFO LS method) occasionally finds reasonable solutions.  

%
%
%

\bibliographystyle{abbrv}
\bibliography{references_KH.bib}

\end{document}